\documentclass[a4paper,11pt,oneside]{article}
\usepackage{amsthm,amsmath,amssymb,amsfonts,amscd}
\usepackage{hyperref}
\hypersetup{
    colorlinks=false,
    linkcolor=blue
}
\usepackage{graphicx}
\usepackage{enumerate}
\usepackage[all]{xy}
\usepackage{booktabs}

\usepackage[english]{babel}

\usepackage[utf8]{inputenc}
\usepackage{amsmath,amsthm,amssymb,caption,amsbsy,graphicx,mathtools,enumerate}

\usepackage[shortlabels]{enumitem}
\usepackage{centernot}
\usepackage{fontenc}
\usepackage{array, multirow}
\usepackage{mdframed,cleveref}
\usepackage{bbm}
\usepackage{tikz,float}
\usepackage[margin=1cm,font=small]{caption}
\usepackage[position=top]{subcaption}

\usepackage{enumerate}
\usepackage{lipsum}
\usepackage{verbatim} 
\usepackage{float}
\usepackage{circuitikz}
\usepackage{caption}
\usepackage{media9}
\usepackage[nottoc]{tocbibind}

\usepackage{geometry} 
\usepackage{blkarray}

\usepackage[symbol]{footmisc}

\usepackage{animate}

\usepackage{wasysym}
\usepackage{fontawesome}

\usepackage{cancel}
\usepackage{ulem} 

\usepackage{multicol}

\usepackage{multirow}

\usepackage{varwidth}

\usepackage{stackengine}
\usepackage{scalerel}

\usepackage{extarrows} 

\usepackage{listings} 
\usepackage{color}

\usepackage{tcolorbox}

\usepackage{nicefrac}

\definecolor{jet1}{rgb}{0, 0, 0.6667}
\definecolor{jet2}{rgb}{0, 0, 1}
\definecolor{jet3}{rgb}{0, 0.3333, 1}
\definecolor{jet4}{rgb}{0, 0.6667, 1}
\definecolor{jet5}{rgb}{0, 1, 1}
\definecolor{jet6}{rgb}{0.3333, 1, 0.6667}
\definecolor{jet7}{rgb}{0.6667, 1, 0.3333}
\definecolor{jet8}{rgb}{1, 1, 0}
\definecolor{jet9}{rgb}{1, 0.6667, 0}
\definecolor{jet10}{rgb}{1, 0.3333, 0}

\usepackage{pgfplots}
\usepackage{tkz-euclide}
\pgfplotsset{compat=1.16}
\pgfplotsset{colormap/jet
}
\usetikzlibrary{intersections}

\newtheorem{theorem}{Theorem}[section]

\newtheorem{remark}[theorem]{Remark}

\usepackage[colorinlistoftodos]{todonotes}

\newcommand{\M}{\mathcal{M}}

\usepackage{comment}
\usepackage{soul}

\title{Traveling waves in a model for cortical spreading depolarization with slow-fast dynamics}
\date{}
\author{David Reyner-Parra$^{1}$, Carles Bonet$^{1}$, Teresa M. Seara$^{1,2,3}$ and Gemma Huguet$^{1,2,3}$ \\
\parbox{12.5cm}{
  \small
  \begin{itemize}
  \item[$^1$]
    Departament de Matem\`atiques, Universitat Polit\`ecnica de Catalunya, Barcelona, Spain 
  \item[$^2$]
    Institut de Matem\`atiques de la UPC - Barcelona Tech (IMTech), Barcelona, Spain 
  \item[$^ 3$]
  Centre de Recerca Matem\`atica, Barcelona, Spain
   \end{itemize}
 }}

\begin{document}

\maketitle

\textbf{Corresponding author:} Gemma Huguet, \texttt{gemma.huguet@upc.edu} \\

\abstract{
Cortical spreading depression and spreading depolarization (CSD) are waves of neuronal depolarization that spread across the cortex, leading to a temporary saturation of brain activity. They are associated to various brain disorders such as migraine and ischemia. We consider a reduced version of a biophysical model of a neuron-astrocyte network for the initiation and propagation of CSD waves (Huguet et al., Biophys. J. , 2016), consisting of reaction-diffusion equations. The reduced model considers only the dynamics of the neuronal and astrocytic membrane potentials and the extracellular potassium concentration, capturing the instigation process implicated in such waves. We present a computational and mathematical framework based on the parameterization method and singular perturbation theory to provide semi-analytical results on the existence of a  wave solution and to compute it jointly with its velocity of propagation. The traveling wave solution can be seen as an heteroclinic connection of an associated system of ordinary differential equations with a slow-fast dynamics.  The presence of distinct time scales within the system introduces numerical instabilities, which we successfully address through the identification of significant invariant manifolds and the implementation of the parameterization method. 
Our results provide a methodology that allows to identify efficiently and accurately the mechanisms responsible for the initiation of these waves and the wave propagation velocity.}

\vspace{0.5cm}


\section{Introduction}
    
Cortical spreading depression and depolarization (CSD) are slowly propagating waves of rapid, near-complete depolarization of neurons and glial cells across the cortex \cite{Dreier11, PietrobonM14, Kramer16}. These waves move slowly through the cortex (2 to 6 mm/min), leading to a temporary cessation of brain electrical activity. This sustained depolarization of neurons has been implicated in various neurological disorders such as migraine and ischemia \cite{Lauritzenetal10, Kramer16}.

The hallmark of CSD is the disruption of transmembrane ionic gradients, which results in the collapse of ion homeostasis, swelling of cells, and elevated extracellular potassium levels \cite{Somjen01, Dreier11, Kramer16}. Astrocytes may play a crucial role in regulating potassium concentration by absorbing its accumulation from the extracellular space \cite{BellotS17}, thus potentially delaying or even preventing the onset of CSD waves \cite{Ma15, ZhaoR10, Nedergaard05}. The ability of astrocytes to regulate potassium levels highlights the importance of astrocyte-neuron interaction in the control of brain function.

    Dynamics involved in these phenomena occur at different time scales. The processes governing ionic concentration changes occur at a slower time scale (seconds) compared to those that activate ionic currents and alter cell membrane potential (milliseconds). 

   In this work, we examine a biophysical model for cortical depolarization developed in \cite{Huguet16}. The model consists of a network of neurons and astrocytes, which was one of the first to incorporate a detailed description of astrocyte functioning. The interaction between neurons and astrocytes occurs through the diffusion of sodium and potassium ions in the extracellular space, as well as the release and uptake of these ions by the cells. This model was used in \cite{Huguet16} to investigate, by means of numerical simulations, how the specific properties of astrocyte processes and the coupling between neurons and astrocytes, interact to generate or prevent the propagation of depolarization waves. 
   
In this paper, we provide a reduction of the model in \cite{Huguet16} based on the methodology developed in \cite{Lee17}. The model is simplified so that it includes only those processes needed to initiate the wave. Thus, the reduced model considers only the dynamics of the neuronal and astrocytic membrane potentials and the extracellular potassium, while  retaining the neuronal and glial pump currents, potassium release during firing, and diffusion along the extracellular space. This simplified framework makes the model more amenable to theoretical analysis.

 We perform a study of the reduced model using computational and analytical tools to provide semi-analytical results on the existence of a traveling wave solution and its propagation velocity. Our approach consists in looking for a traveling wave solution by searching for a heteroclinic orbit in the resulting system of ordinary differential equations. The presence of different times scales in the system  allows us to apply slow-fast theory to study the system. Thus, we adopt the singular perturbation setting \cite{Kuehn15, Jones_book} and introduce a small parameter that represents the ratio between the different time scales. Using techniques from classical singular perturbation theory \cite{Kuehn15, Jones_book}, we break down the system into its fast and slow subsystems and construct an heteroclinic orbit in the singular limit. By leveraging the geometric understanding of the dynamics in the singular limit, we design a numerical strategy to compute the heteroclinic orbit beyond this limit, which effectively overcomes the numerical instabilities arising from the inherent time-scale differences in the system. In our study, we present two distinct numerical strategies to achieve this, one based on the parameterization method \cite{CabreFL05, Haro_book} and another employing Fenichel's theory, which ensures the persistence of normally hyperbolic invariant manifolds \cite{Fenichel71,Fenichel73, HPS77}. 

 Both methods have proven successful in mitigating the numerical instabilities associated with the system's time-scale differences. They yield accurate results for computing the traveling wave and its propagation velocity. Importantly, our methods shed light on the underlying mechanisms responsible for initiating the wave, providing valuable insights that can potentially be extended to the full model in the future.

The paper is organized as follows. In Section~\ref{sec:model} we present the original model introduced in \cite{Huguet16} to describe CSD and the reduced version capturing the initiation of the wave. In Section~\ref{sec:analysis} we present the mathematical analysis of the traveling wave solutions of the reduced model using the singular perturbation framework and we construct a singular heteroclinic orbit. In Section~\ref{sec:full_system} we provide the computation of the heteroclinic orbit and the velocity  for the reduced system using the parameterization method \cite{Haro_book, CabreFL05}. In 
    Section~\ref{sec:fenichel} we provide an alternative method to compute the heteroclinic orbit numerically based on Fenichel's theory \cite{Fenichel71, Fenichel73}.  We end with a discussion in Section~\ref{sec:discussion}. Finally, the Appendix contains some technical details involving algebraic mathematical computations.

    \section{The mathematical model}\label{sec:model}

    \subsection{Model for cortical spreading depression involving astrocytes}

    In this section, we present the neuron-astrocyte network model introduced in  \cite{Huguet16}, upon which we base our work. The model describes quantitatively the dynamics of a network consisting of one array of neurons and another of astrocytes, sharing only the extracellular medium, through which there is diffusion of sodium ($Na^+$) and potassium ($K^+$) ions. The model captures accurately the initiation and propagation of cortical spreading depolarization (CSD) waves, emerging in several cortical diseases, such as migraine and ischemia \cite{Dreier11}. A thorough study on preventing or delaying this phenomenon was carried out in \cite{Huguet16}, with particular emphasis on the role played by the astrocyte cells and their gap junctions.

    The mathematical model consists of a set of ordinary and partial differential equations modeling the dynamics of the neuronal and astrocytic membrane potentials, $V_N$ and $V_A$, respectively, the extracellular and intracellular concentrations of sodium and potassium, $[Na^+]_e$, $[Na^+]_i$, $[Na^+]^A_i$, $[K^+]_e$, $[K^+]_i$, $[K^+]^A_e$, respectively, and the gating variables $n$ and $h_p$. The system of differential equations is given by:
    
    \begin{equation} \label{gHuguet_model}
        \begin{split}
            C_m \frac{d V_N}{d t} &= -I_{Na} - I_{NaP} - I_K - I_L - I_{Pm} \, ,\\[0.2em]
            C_m^A \frac{d V_A}{d t} &= -I_{Na}^A - I_K^A - I_{Pm}^A - I_{\text{gap}} \, ,\\[0.2em]
            \frac{d n}{d t} &= \phi_n \frac{n_\infty(V_N) - n}{\tau_n (V_N)} \, , \\[0.2em]
            \frac{d h_p}{d t} &= \phi_{h_p} \frac{{h_p}_\infty(V_N) - h_p}{\tau_{h_p}(V_N)} \, , \\[0.2em]
            \frac{d [Na^+]_i}{d t} &= -\frac{10 S_N}{F \Omega_n} (I_{Na} + I_{NaP} + 3 I_{Pm}) \, , \\[0.2em]
            \frac{d [Na^+]_i^A}{d t} &= -\frac{10 S_A}{F \Omega_a} (I_{Na}^A + 3 I_{Pm}^A - I_{Na,\text{gap}}) \, , \\[0.2em]
            \frac{d [K^+]_i}{d t} &= -\frac{10 S_N}{F \Omega_n} (I_K - 2 I_{Pm}) \, , \\[0.2em]
            \frac{d [K^+]_i^A}{d t} &= -\frac{10 S_A}{F \Omega_a} (I_K^A - 2 I_{Pm}^A - I_{K,\text{gap}}) \, , \\[0.2em]
            \frac{\partial [Na^+]_e}{\partial t} &= D_{Na} \frac{\partial^2 [Na^+]_e}{\partial x^2} + \frac{10 S_N}{F \Omega_e} (I_{Na} + I_{NaP} + 3 I_{Pm}) + \frac{10 S_A}{F \Omega_e} (I_{Na}^A + 3 I_{Pm}^A) \, , \\[0.2em]
            \frac{\partial [K^+]_e}{\partial t} &= D_K \frac{\partial^2 [K^+]_e}{\partial x^2} + \frac{10 S_N}{F \Omega_e} (I_K - 2 I_{Pm}) + \frac{10 S_A}{F \Omega_e} (I_K^A - 2 I_{Pm}^A) \, . \\[0.2em]
        \end{split}
    \end{equation}
    
    The ionic currents $I_{Na}$, $I_{NaP}$, $I_{K}$ and $I_L$ are modelled as in the Hodgkin-Huxley model, that is,
    \begin{equation}
        \begin{split}
            I_{Na} &= g_{Na} \, m_\infty^3(V_N) \, (1 - n) (V_N - E_{Na}) \, , \\[0.2em]
            I_{NaP} &= g_{NaP} \, m_{p \infty}^3(V_N) \, h_p (V_N - E_{Na}) \, , \\[0.2em]
            I_{K} &= g_{K} \, n^4 (V_N - E_{K}) \, , \\[0.2em]
            I_L &= g_L (V_N - E_L) \, .
        \end{split}
    \end{equation}
    The reversal potential $E_X$, for $X \in \{Na, K\}$, is 
    modelled by the Nernst equation
    \begin{equation} \label{nernst_potential}
        E_X = \frac{R T}{F} \text{ln} \frac{[X]_e}{[X]_i}, \, 
    \end{equation}
where $R,T$, and $F$ are the gas constant, temperature, and Faraday's constant, respectively. The values are given in Table~\ref{tab:parameters}. 
    
    The asymptotic functions $m_\infty$, $m_{p\infty}$, $n_\infty$ and $h_{p \infty}$ modeling the steady state values of the ionic gating variables are sigmoidal-shaped functions of the form
    \begin{equation}
        X_\infty(V) = \frac{1}{1 + e^{(-(V - V_X)/\theta_X)}} \, ,
    \end{equation}
    for $X \in \{m,n,m_p,h_p \}$, with $V_{m} = -34$, $\theta_m = 5$, $V_{n} = -55$, $\theta_n = 14$, $V_{m_p} = -40$, $\theta_{m_p} = 6$, $V_{h_p} = -48$ and $\theta_{h_p} = -6$. Finally, the time constants modeling the rate at which the gating variables reach its steady state values are
    \begin{equation}
        \begin{split}
            \tau_{n}(V) &= 0.05 + \frac{0.27}{1 + e^{(-(V + 40)/(-12))}} \, , \\[0.2em]
            \tau_{h_p}(V) &= \frac{10000}{\text{cosh}((V + 48)/12)} \, .
        \end{split}
    \end{equation}

    The ionic currents $I_{Na}^A$ and $I_K^A$ in the modelization of the astrocyte equation are described by the Goldman-Hodgkin-Katz (GHK) formulation:
    \begin{equation}\label{eq:astro_currents}
        \begin{split}
            I_{Na}^A &= P_{Na} F \phi \frac{[Na^+]_e \, e^{-\phi} - [Na^+]_i^A}{e^{-\phi} - 1} \, , \\[0.2em]
            I_{K}^A &= P_K F \phi \frac{[K^+]_e \, e^{-\phi} - [K^+]_i^A}{e^{-\phi} - 1} \, ,
        \end{split}
    \end{equation}
    with $\phi := \frac{F}{R T}V_A$, $P_{Na} = 1.5 \cdot 10^{-8}$ and $P_K = 1 \cdot 10^{-6}$. Notice that the currents $I^A_{Na}$ and $I^A_K$ in \eqref{eq:astro_currents} are not defined when the astrocyte's membrane potential reaches zero because the dividing term $e^{-\phi} - 1$ cancels. This is resolved by explicitly defining the values of $I_{Na}^A$ and $I_{K}^A$ at $V_A = 0$,
    \begin{equation*}
        \begin{split}
            I_{Na}^A(V_A,[Na^+]_e, [Na^+]_i^A)_{|_{V_A = 0}} & := \lim_{V_A \rightarrow \, 0} \, P_{Na} F \phi \frac{[Na^+]_e \, e^{-\phi} - [Na^+]_i^A}{e^{-\phi} - 1} = P_{Na} F \big([Na^+]_i^A - [Na^+]_e\big) \, , \\[0.2em]
            I_K^A(V_A,[K^+]_e, [K^+]_i^A)_{|_{V_A = 0}} & := \lim_{V_A \rightarrow \, 0} \, P_K F \phi \frac{[K^+]_e \, e^{-\phi} - [K^+]_i^A}{e^{-\phi} - 1} = P_{K} F \big([K^+]_i^A - [K^+]_e\big) \, .
        \end{split}
    \end{equation*}

    The currents $I_{Na,\text{gap}}$ and $I_{K,\text{gap}}$ in the astrocytes correspond to the potassium and sodium currents exchanged with nearby gap-junction-coupled astrocytes. They are modelled using the GHK formalism in \cite{Huguet16}. Since in this paper we will not consider them, we do not include details of the modelization here and we refer the reader to \cite{Huguet16} for further information.

    The sodium-potassium pump currents in neurons and astrocytes, $I_{Pm}$ and $I_{Pm}^A$, respectively, are given by
    \begin{equation}
        \begin{split}
            I_{Pm} &= \rho_{N} \bigg( \frac{[K^+]_e}{2+ [K^+]_e} \bigg)^2 \bigg( \frac{[Na^+]_i}{7.7 + [Na^+]_i} \bigg)^3 \, , \\[0.2em]
            I_{Pm}^{A} &= \rho_{A} \bigg( \frac{[K^+]_e}{2 + [K^+]_e} \bigg)^2 \bigg( \frac{[Na^+]_i^{A}}{7.7 + [Na^+]_i^{A}} \bigg)^3 \, ,
        \end{split}
    \end{equation}
     where $\rho_N, \rho_A$ represent the maximal current allowed to go through the neuron and astrocyte pumps, respectively (regarded as parameters in \cite{Huguet16} to explore different conditions). In the present work we have fixed both of them at a constant value of 5, for which we know from \cite{Huguet16} that a wave can be initiated and propagated through the tissue. We know from previous works \cite{Huguet16, Lee17} that increasing the pump strength, $\rho_N$ and $\rho_A$, increases the $[K^+]_e$ threshold for wave initiation, which yields either a delay in the wave initiation (it takes longer to achieve the threshold) or a propagation failure (when the threshold is not achieved).

\begin{figure}[htbp!]
        \centering
        \includegraphics[width=\textwidth]{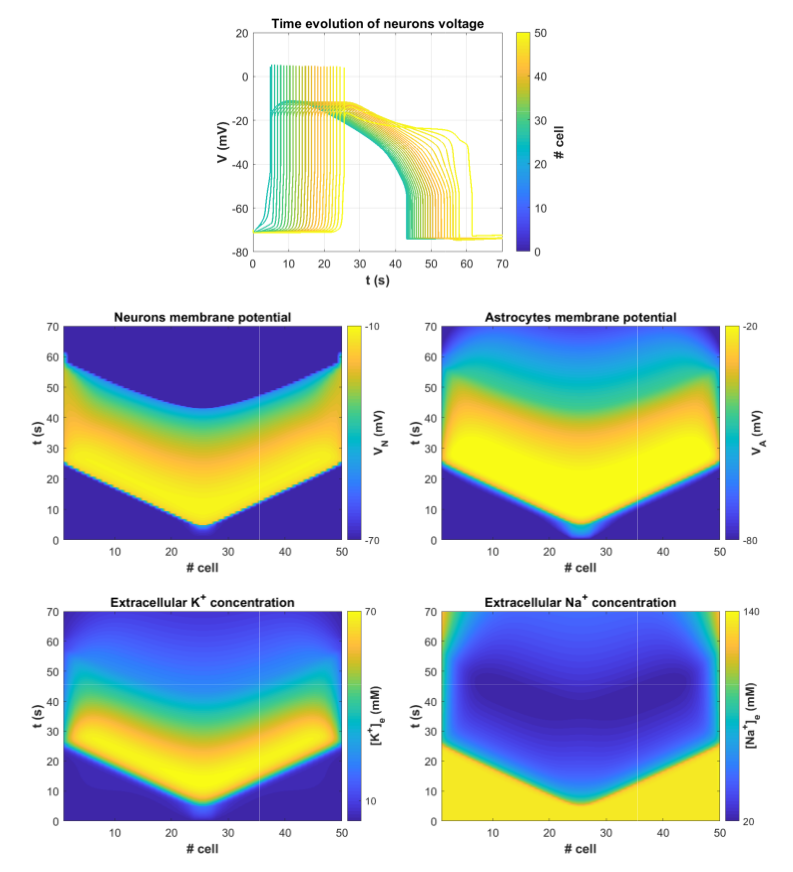}
        \caption{Solutions of the model \eqref{gHuguet_model} with simulated $K^+$ injection. The simulation consists of a 1D array network of 50 neuron-astrocyte pairs sharing the same extracellular space, with the middle four cells externally stimulated with an insult of $K^+$ at a rate of $5$ mM/s until a depolarizing wave is triggered. The initial conditions are identical for all the pairs. The diffusive terms in system \eqref{gHuguet_model} are approximated by central second-order finite differences, imposing Dirichlet boundary conditions on the endpoints of the array (using typical levels of $K^+$ and $Na^+$ in healthy tissues). (Top) Temporal profile of the variable $V_N$ for each cell during a CSD wave (each curve is coloured according to the neuron index). (Bottom) 2D views of the traveling wave solution (arising) on $V_N$, $V_A$, $[K^+]_e$ and $[Na^+]_e$, spreading through space from middle cells towards adjacent ones over time. Parameters are given in Table \ref{tab:parameters}.}
        \label{simulations_ghuguet_model}
    \end{figure}

    We performed a simulation of system \eqref{gHuguet_model} with the parameters given in Table~\ref{tab:parameters} for a network of 50 neuron-astrocyte pairs sharing the extracellular space.
    To do so, we discretize the space variable $x$ using a grid of points $x_i = i \Delta x$ for $i = 1,\dots,N$, where $\Delta x$ is the distance between neurons. By employing finite differences and imposing Dirichlet boundary conditions  at the endpoints of the array (which aims to simulate healthy tissue), we obtain a system of ODEs. See Appendix~\ref{ap:pde} for more details. We start the network at rest and inject an insult of $K^+$ to the extracellular space of the middle cells until a depolarizing wave is evoked. More precisely, the insult is modeled by  adding a positive constant, $I_{K,\text{inj}}=0.005$ to the right-hand side of the equation for $[K^+]_e$ corresponding to the middle cells until their voltage reaches the value $-30$mV.   The simulation shows a depolarization wave which propagates throughout the array of neuron-astrocyte pairs (see Figure~\ref{simulations_ghuguet_model}).  
    
    This solution corresponds to a traveling wave, a special type of solution of a partial differential equation (PDE) on an infinite domain with a strong restraint between $\boldsymbol{x}$ and $\boldsymbol{t}$ (i.e. $\boldsymbol{x}\pm c \boldsymbol{t}$), whose spatial profile advances in time at a constant speed $c$ (therefore being invariant by temporal translations).

    \subsection{Model reduction}

    In this section, we present a reduction of model \eqref{gHuguet_model} following the strategy employed in \cite{Lee17}. First, we provide details of the simplifications and its main differences with respect to \cite{Lee17}. Then, we present simulations of the reduced model, showing that it retains the ability to initiate and propagate a depolarizing front along the array of neuron-astrocyte pairs.

    The strategy for the model reduction is to retain only the dynamics for the variables that play a significant role in the initiation of CSD waves, while neglecting the dynamics of the rest of variables. Thus, since the variation of sodium concentration is not essential for the wave initiation and propagation,  the variables  $[Na^+]_e$, $[Na^+]_i^A$ and $[Na^+]_i$ are frozen at the resting values.  
    Namely, the resting values are set at 3.5 (mM) for the intracellular concentrations (i.e., $[Na^+]_i=[Na^+]_i^A =3.5$) and at 135 (mM) for the extracellular one (i.e. $[Na^+]_e=135$).

    On the other hand, since the astrocytes' volume is bigger than that of the extracellular medium (i.e. $\Omega_e = \alpha_0 (\Omega_n + \Omega_a)$, $\alpha_0 \ll 1$), changes in the quantity of  $K^+$ ions will have a greater impact on the extracellular $K^+$ concentration $[K^+]_e$ than on the concentration within the astrocytes $[K^+]_i^A$. Thus, the variable $[K^+]_i^A$ in equation \eqref{gHuguet_model} is also taken as constant. Finally, the conservation of $K^+$ ions states that
    \begin{equation*}
        [K^+]_i := \frac{1}{\Omega_n} ( K^+_{\text{tot}} - \Omega_e [K^+]_e - \Omega_a [K^+]_i^A) \, ,
    \end{equation*}
    being $K^+_{\text{tot}}$ the total initial number of $K^+$ ions in the space. Note that, unlike $[Na^+]_i$, the magnitude $[K^+]_i$ is not constant as it depends moderately on the slowly changing dynamics of $[K^+]_e$. The initial values for the extracellular and intracellular $K^+$ concentrations are $[K^+]_e=3.5$ (mM) and $[K^+]_i=135$ (mM), respectively, while $[K^+]^A_i$ is kept constant at $135$ (mM). 

    The model is further simplified by removing the dynamics of the two gating variables $n$ and $h_p$. The variable $n$ has a small time constant, so it approaches fast the steady state $n_{\infty}$, thus $n$ is set to $n_\infty$. On the other hand, $h_p$ is a very slow variable, so we freeze it at its resting value. As we shall see later, the dynamics of $h_p$ plays a role in restoring the polarization of the neuron and astrocyte cells. When $h_p$ is held constant, the reduced system will show depolarizing fronts instead of pulses, which is sufficient to study the initiation and propagation of the waves.

    We know from \cite{Huguet16, Ma15} that the astrocytes form a syncytium through gap junctions, which have a critical role in preventing wave propagation. Indeed, when the strength of the gap junctions is weak, a wave is initiated, as depicted in Figure \ref{simulations_ghuguet_model}. 
    However, in this paper, our focus is not on incorporating the effects of gap junctions. Instead, we intentionally remove the current through the astrocyte gap junctions $I_{\text{gap}}$ from the original model. This  allows us to explore mathematical algorithms that capture the properties of the initiated wave.
    
    These simplifications render system \eqref{gHuguet_model} depending only on the variables $V_N$, $V_A$ and $[K^+]_e$ in the following way:
    \begin{equation} \label{current_pde_model}
        \begin{split}
            C_m \frac{d V_N}{d t} &= -I_{Na} - I_{NaP} - I_K - I_L - I_{Pm} \, ,\\[0.2em]
            C_m^A \frac{d V_A}{d t} &= -I_{Na}^A - I_K^A - I_{Pm}^A  
            \, ,\\[0.2em]
            \frac{\partial [K^+]_e}{\partial t} &= D_K \frac{\partial^2 [K^+]_e}{\partial x^2} + \frac{10 S_N}{F \Omega_e} (I_K - 2 I_{Pm}) + \frac{10 S_A}{F \Omega_e} (I_K^A - 2 I_{Pm}^A) \, . \\[0.2em]
        \end{split}
    \end{equation}
    System \eqref{current_pde_model} is the model we will analyze in this paper, which is the same considered in \cite{Lee17} except that we keep the terms $I_K^A$ and $I_{Pm}^A$ in the equation for $[K^+]_e$, describing the influence of astrocytes in the extracellular levels of $K^+$.

\begin{figure}[htbp!]
         \includegraphics[width=\textwidth]{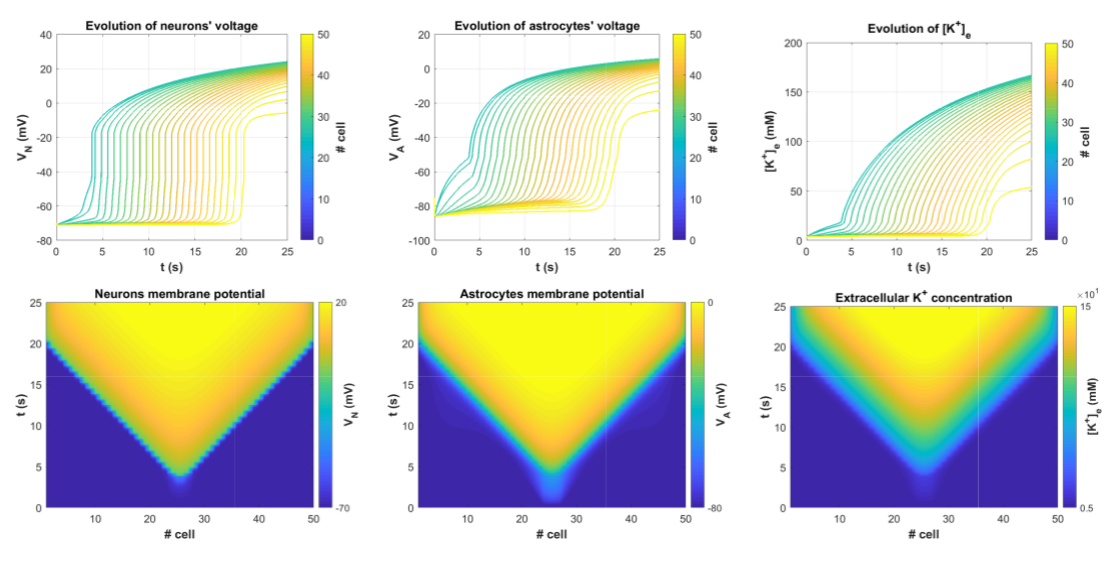}
        \caption{Solutions of model \eqref{current_pde_model} with simulated $K^+$ injection. Fifty neuron-astrocyte pairs are simulated, with the middle four cells stimulated
with insult of $K^+$ at a rate of $5$ mM/s until a front of spreading depolarization is evoked. (Top) Temporal profile of the variables $V_N$, $V_A$ and $[K^+]_e$ for each cell during a front of spreading depolarization (each curve is coloured according to the neuron-astrocyte index). (Bottom) 2D views of a depolarizing front on $V_N$, $V_A$ and $[K^+]_e$ spreading through space to each neighbouring cell over time. As in Figure \ref{simulations_ghuguet_model}, Dirichlet boundary conditions have been applied on the extracellular potassium $K^+$, to simulate typical levels of $[K^+]_e$ in healthy tissues. Parameters are given in Table \ref{tab:parameters}.}
        \label{simulations_rlee_reduced_gh_model}
    \end{figure}

    In Figure \ref{simulations_rlee_reduced_gh_model} we show a numerical simulation of the reduced model \eqref{current_pde_model}. As mentioned before, the resulting system can capture the initiation of the depolarization wave but fails to capture the cell's membrane potential recovery phase to the resting potential, as the variable $h_p$ is frozen. Thus, it does no longer admits traveling pulse solutions (as displayed in Figure \ref{simulations_ghuguet_model}) but traveling fronts. As in simulations of Figure \ref{simulations_ghuguet_model}, we modelled an insult of $K^+$ in the extracellular space around the four middle neurons, in order to initiate the front. \\

\begin{figure}[htbp!]
        \centering
        \includegraphics[width=0.6\textwidth]{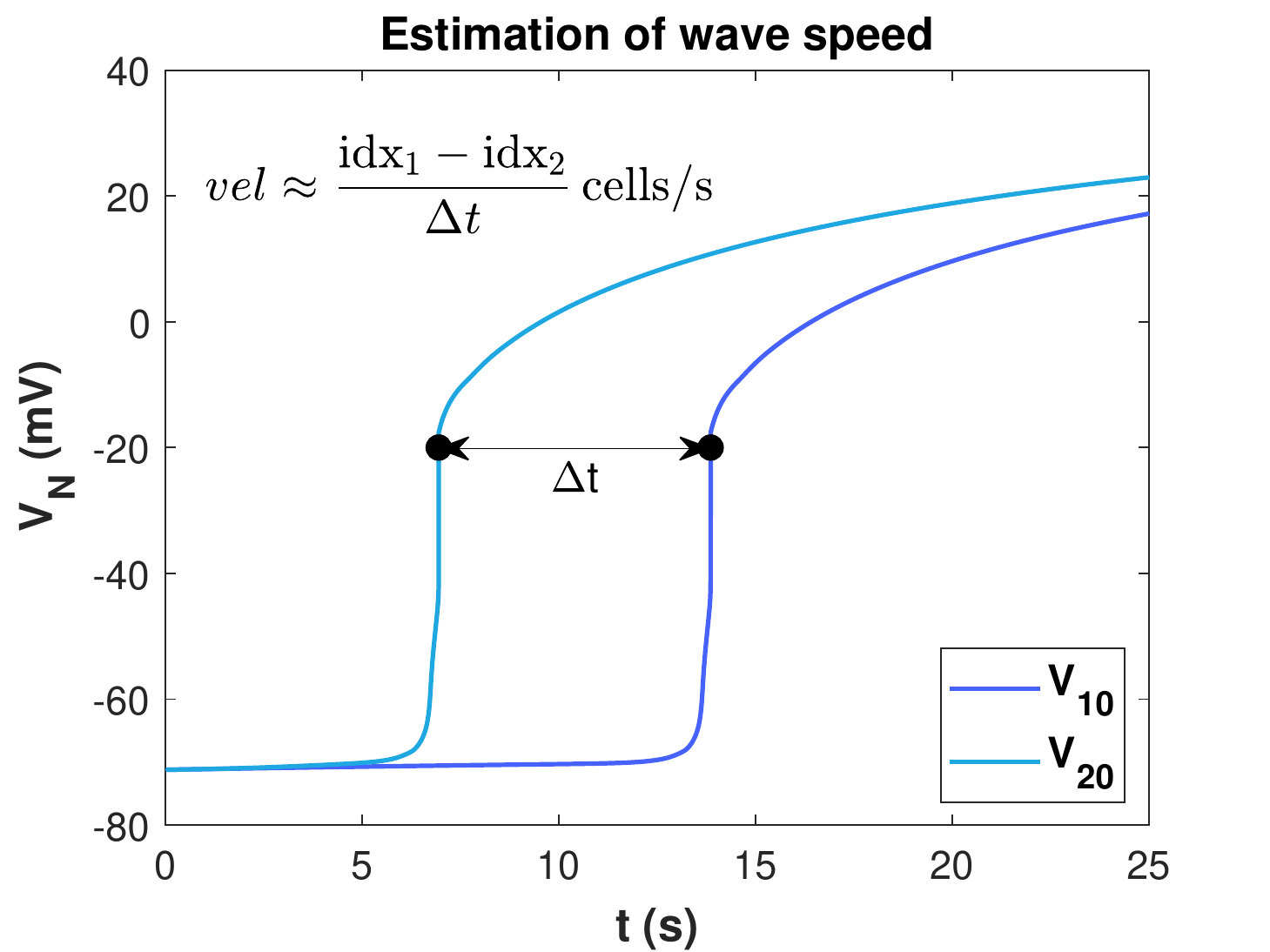}
        \caption{Estimation of the wave speed of the numerical simulations of system \eqref{current_pde_model} shown in Figure \ref{simulations_rlee_reduced_gh_model} from the voltage traces of neurons 10 and 20. Using that the distance between cells is $\Delta x= 0.044 \, \text{mm}$, the estimated front velocity is $vel \approx 3.82 \, \text{mm}/\text{min}$.}
        \label{wave_speed_estimation}
    \end{figure}

    From the numerical simulations, we can provide an estimate of the wave speed of the traveling front of system \eqref{current_pde_model} by computing the elapsed time between the depolarization of two distinct cells. The speed is then obtained as the distance between the two cells divided by this time (see Figure \ref{wave_speed_estimation}). We want to emphasize that we can rely on the numerical simulations to establish the existence of a traveling wave solution, but accurate estimates of its speed cannot be made due to the influence of the boundary conditions. Although we estimate the velocity using the voltage of cells that are far from the boundaries, we can still observe that the chosen boundary conditions (modeling healthy tissue at the boundaries of the domain) impact on the propagation speed of the traveling front by slowing it down. Indeed, our numerical explorations show that enlarging the size of the medium of the simulation by considering more neuron-astrocyte pairs, increases significantly the average wave speed (see Table \ref{tab:table_wave_speeds}).  

    \begin{table}[htbp!]
        \centering
        \begin{tabular}{|c|c|}
            \hline
            \textbf{\# neuron-astrocyte pairs} & \textbf{Wave speed estimation (mm/min)} \\
            \hline
            50 & 3.8205 \\
            100 & 4.1641 \\
            300 & 4.7186 \\
            500 & 4.8253 \\
            \hline
        \end{tabular}
        \caption{Estimated wave speed of numerical simulations of system \eqref{current_pde_model} as the number of neuron-astrocyte pairs in the medium increases.}
        \label{tab:table_wave_speeds}
    \end{table}

    A further simplification of model \eqref{current_pde_model} has been done in \cite{Lee17}, by leveraging the different scales in the rate of change of the membrane potential variables, $V_N$ and $V_A$, and the extracellular $K^+$ (compare, for instance, the evolution of $V_N$ and $V_A$ relative to $[K^+]_e$ in Figure \ref{simulations_rlee_reduced_gh_model}). Mathematically, imposing that
    \begin{equation}\label{eq:Lee_red}
        \begin{split}
            -I_{Na} - I_{NaP} - I_K - I_L - I_{Pm} &= 0 \, ,\\[0.2em]
            -I_{Na}^A - I_K^A - I_{Pm}^A &= 0 \, ,\\[0.2em]
        \end{split}
    \end{equation}
    the system is reduced to a single PDE for $[K^+]_e$ of reaction-diffusion type, that preserves qualitatively the dynamics shown in Figure \eqref{simulations_rlee_reduced_gh_model}. As we will see later, part of our study of system \eqref{current_pde_model} will require the knowledge of the dynamics of this further reduced model. 
    
    \begin{table}
    \begin{center}
        \begin{tabular}{|c |c | c | p{8cm}|}
         
            \hline
            Parameters & Paper & Units & Description 
            \\
            \hline
            R & 8.31 & J /mol K & ideal gas constant 
            \\
            F & 96485 & C/mol & Faraday's constant \\
            T & 310 & K & absolute temperature 
            \\
            $C_m$ & 1 & $\mu$F/cm$^2$ & neuron's membrane capacitance per unit area 
            \\
            $C^A_m$ & 1 & $\mu$F/cm$^2$ & astrocyte's membrane capacitance per unit area 
            \\
            $\phi_n$ & 0.8 & ms$^{-1}$ & maximum rate of activation of $K^+$ channels 
            \\
            $\phi_{h_p}$ & 0.05 & ms$^{-1}$ & maximum rate of activation of $Na^+$ channels 
            \\
            $g_{Na}$ & 50.0 & mS/cm$^2$ & maximum fast $Na^+$ current conductance 
            \\
            $g_{NaP}$ & 0.8 & mS/cm$^2$ & maximum fast $Na^+$ current conductance 
            \\
            $g_K$ & 15 & mS/cm$^2$ & maximum delayed rectifier $K^+$ conductance 
            \\
            $g_L$ & 0.5 & mS/cm$^2$ & nonspecific leak current conductance 
            \\
            $E_L$ & -70 & mV & leak current reversal potential 
            \\
            $S_N$ & 922 & $\mu$m$^2$ & neuron membrane total surface area 
            \\
            $\Omega_n$ & 2160 & $\mu$m$^3$ & neuron total volume 
            \\
            $S_A$ & 1600 & $\mu$m$^2$ & astrocyte membrane total surface area 
            \\
            $\Omega_a$ & 2000 & $\mu$m$^3$ & astrocyte total volume 
            \\
            $\Omega_e$ & $ \alpha_0 (\Omega_n + \Omega_a)$ & $\mu$m$^3$ & extracellular volume 
            \\
            $\alpha_0$ & 0.2 & & volume fraction 
            \\
            $D_K$ & $1.96 \cdot 10^{-5}$ & $cm^2/s $ & $K^+$ diffusion coefficient 
            \\
            $\Delta x$ & 0.044 &  mm & distance between neurons 
            \\
            $D_{Na}$ & $1.33 \cdot 10^{-5}$ &  cm$^2$/s & Na$^+$ diffusion coefficient \\
            $P_K$ & $1 \cdot 10^{-6}$ & cm/s & membrane K$^+$ permeability coefficient \\
            $P_{Na}$ & $1.5 \cdot 10^{-8}$ & cm/s & membrane Na$^+$ permeability coefficient 
            \\
            $\rho_N$ & 5 & $\mu$A/cm$^2$ & maximum neuron Na$^+$-K$^+$ ATPase pump current \\
            $\rho_A$ & 5 & $\mu$A/cm$^2$ & maximum astrocyte Na$^+$-K$^+$ ATPase pump current 
            \\
            \hline
        \end{tabular}
     \end{center}
     \caption{Parameters for system \eqref{gHuguet_model} used in this paper which combine parameter values from \cite{Huguet16} and \cite{Lee17}. Note that $E_X$ in \eqref{nernst_potential} must be multiplied by a factor $1000$ to convert from V to mV.} 
     \label{tab:parameters}
    \end{table}

    \section{Traveling wave solutions in the singular limit} \label{sec:analysis}
    
    In this section we present a strategy that combines analytical and numerical methods to show that system \eqref{current_pde_model} has 
    traveling wave solutions and to compute them. 
    The first step is to postulate a solution of the form 
    \begin{equation}\label{eq:ansatz}
        \begin{split}
            V_N(\boldsymbol{x}, t) &= \overline{V}_N(\xi) \, , \\[0.2em]
            V_A(\boldsymbol{x}, t) &= \overline{V}_A(\xi) \, , \qquad \text{with} \quad \, \xi := \boldsymbol{x}/\sqrt{D_K} + c \, t \\[0.2em]
            [K^+]_e(\boldsymbol{x}, t) &= \overline{[K^+]}_e(\xi) \, ,
        \end{split}
    \end{equation}
    where the constant $c$ (to be determined) establishes the wave propagation velocity, $vel=c \sqrt{D_K}$. Substituting \eqref{eq:ansatz} in \eqref{current_pde_model} results in a system of ordinary differential equations, given by 
    \begin{equation} \label{neuron_astrocyte_model_wave_sol}
        \begin{split}
            c \frac{d \overline{V}_N}{d \xi} &= -\frac{1}{C_m} (I_{Na} + I_{NaP} + I_K + I_L + I_{Pm}) \, , \\[0.2em]
            c \frac{d \overline{V}_A}{d \xi} &= -\frac{1}{C_m^A} (I_{Na}^A + I_{K}^A + I_{Pm}^A) \, , \\[0.2em]
            c \frac{d \overline{[K^+]}_e}{d \xi} &= 
            \frac{d^2 \overline{[K^+]}_e}{d \xi^2} + \frac{10 S_N}{F \Omega_e} (I_K - 2 I_{Pm}) + \frac{10 S_A}{F \Omega_e} (I_{K}^A - 2 I_{Pm}^A) \, .
        \end{split}
    \end{equation}
    We will see that system \eqref{neuron_astrocyte_model_wave_sol} has two saddle points $p_{l_1}$ and $p_r$ (see Table \ref{tab:eq_points}). To find a traveling front solution of system \eqref{current_pde_model}, we will look for an heteroclinic orbit of system \eqref{neuron_astrocyte_model_wave_sol} from $p_{l_1}$ to $p_r$. The goal is to determine whether the parameter $c$ can be chosen so that this heteroclinic orbit exists. We will do the analysis for $c>0$. 
    
    We relabel the variables $\overline{V}_N$, $\overline{V}_A$ and $\overline{[K^+]}_e$ as $x$ (not to be confused with the spatial variable $\boldsymbol{x}$), $y$ and $z$, respectively, and introduce the variable $w:=\frac{d \overline{[K^+]}_e}{d \xi}=\frac{dz}{d \xi}$, to finally obtain 
    \begin{equation} \label{neuron_astrocyte_model_wave_sol_4D}
        \begin{split}
            c \frac{d x}{d \xi} &= -\frac{1}{C_m} (I_{Na} + I_{NaP} + I_K + I_L + I_{Pm}) \, , \\[0.2em]
            c \frac{d y}{d \xi} &= -\frac{1}{C_m^A} (I_{Na}^A + I_{K}^A + I_{Pm}^A) \, , \\[0.2em]
            \frac{d z}{d \xi} &= w \, , \\[0.2em]
            \frac{d w}{d \xi} &= 
            c w - \frac{10 S_N}{F \Omega_e} (I_K - 2 I_{Pm}) - \frac{10 S_A}{F \Omega_e} (I_{K}^A - 2 I_{Pm}^A) 
            \, ,
        \end{split}
    \end{equation}
    which, for simplicity, can be written in a more general form
    \begin{equation} \label{neuron_astrocyte_model_wave_sol_4D_simple}
        \begin{split}
            \frac{d x}{d \xi} &= \frac{1}{c} \widetilde f(x,z)\, , \\[0.2em]
            \frac{d y}{d \xi} &= \frac{1}{c} \widetilde g(y,z) \, , \\[0.2em]
            \frac{d z}{d \xi} &= w \, , \\[0.2em]
            \frac{d w}{d \xi} &=  c w - h(x,y,z) \, ,
        \end{split}
    \end{equation}
 
    with
    \begin{equation}\label{eq:fgh}
        \begin{split}
            \widetilde f(x,z) &= -\frac{1}{C_m} \big(I_{Na}(x) + I_{NaP}(x) + I_K(x,z) + I_L(x) + I_{Pm}(z) \big) \, , \\[0.2em]
            \widetilde g(y,z) &= -\frac{1}{C_m^A} \big(I_{Na}^A(y) + I_{K}^A(y,z) + I_{Pm}^A(z) \big) \, , \\[0.2em]
            h(x,y,z) &= \frac{10 S_N}{F \Omega_e} (I_K(x,z) - 2 I_{Pm}(z)) + \frac{10 S_A}{F \Omega_e} (I_{K}^A(y,z) - 2 I_{Pm}^A(z)) \, ,
        \end{split}
    \end{equation}
    where we have explicitly indicated the dependence on $x$, $y$ and $z$ for each ionic current. 

    \begin{figure}[p]
        \centering
        \includegraphics[width=\textwidth]{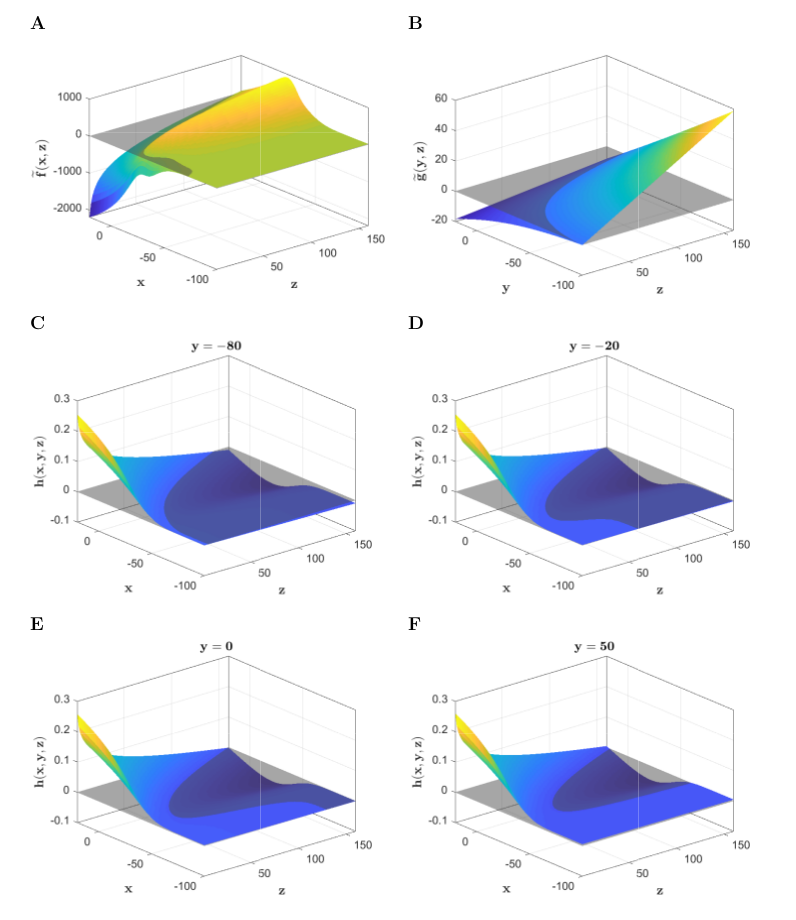}
        \caption{Plots of the graphs of functions given in \eqref{eq:fgh}: (A) $\widetilde f(x,z)$, (B) $\widetilde g(y,z)$ and (C-F) $ h(x,y,z)$ for fixed values of $y$ (as indicated in each panel). The zero-level curve for each function is displayed as the intersection of the grey plane with each surface. Notice also the different orders of magnitude on the value ranges that each function takes, suggesting distinct scales present in system \eqref{neuron_astrocyte_model_wave_sol_4D_simple_eps1}. Parameters of the functions are taken from Table \ref{tab:parameters}.}
        \label{functions_fgh}
    \end{figure}

     We observe that the neuronal and astrocytic membrane potentials ($x$ and $y$ variables, respectively) change more rapidly than that of the extracellular concentration of potassium in the medium ($z$ variable). In Figure \ref{functions_fgh} we plot functions $\widetilde f(x,z)$, $\widetilde g(y,z)$ and $h(x,y,z)$ to illustrate so.  Indeed, the ranges of $\widetilde f$ and $ \widetilde g$ are larger than that of $h$, entailing a higher rate of change in the derivatives $\frac{d x}{d \xi}$ and $\frac{d y}{d \xi}$ and thus a much faster dynamics. For this reason, we  consider system \eqref{neuron_astrocyte_model_wave_sol_4D_simple} with an artificial small parameter $\varepsilon$ to emphasize these slow-fast dynamics, that is,
    \begin{equation} \label{neuron_astrocyte_model_wave_sol_4D_simple_eps}
        \begin{split}
            \mathbf{\varepsilon} \frac{d x}{d \xi} &= \frac{1}{c} f(x,z)\, , \\[0.2em]
            \mathbf{\varepsilon} \frac{d y}{d \xi} &= \frac{1}{c} g(y,z) \, , \\[0.2em]
            \frac{d z}{d \xi} &= w \, , \\[0.2em]
            \frac{d w}{d \xi} &= c w - h(x,y,z)  \, ,
        \end{split}
    \end{equation}
    where we have denoted $f=\varepsilon \widetilde f$ and $g=\varepsilon \widetilde g$. 
    Equivalently, re-scaling the time variable $\varepsilon \eta = \xi$,  system \eqref{neuron_astrocyte_model_wave_sol_4D_simple_eps} becomes
    \begin{equation} \label{neuron_astrocyte_model_wave_sol_4D_simple_eps1}
        \begin{split}
            \frac{d x}{d \eta} &= \frac{1}{c} f(x,z)\, , \\[0.2em]
            \frac{d y}{d \eta} &= \frac{1}{c} g(y,z) \, , \\[0.2em]
            \frac{d z}{d \eta} &= \mathbf{\varepsilon} w \, , \\[0.2em]
            \frac{d w}{d \eta} &= \mathbf{\varepsilon}  ( \, c w - h(x,y,z)  \, ).
        \end{split}
    \end{equation}

\begin{remark}\label{rem:time_scales}
    By looking at the range of values taken by $\widetilde g$ in Figure \ref{functions_fgh}, one can argue that there exists an intermediate time scale between those of $x$ and $z$. This would give rise to another fast-slow version of system \eqref{neuron_astrocyte_model_wave_sol_4D_simple} with two parameters $\varepsilon$ and $\delta$, which would have a term $\varepsilon \, \delta$ multiplying $\dot{x}$ (to indicate the dynamics of $x$ are much faster than that of $y$). In this paper, though, we focus on the existence of two time scales and we leave this analysis for future work.
\end{remark}
    
    Note that equations \eqref{neuron_astrocyte_model_wave_sol_4D_simple_eps} and \eqref{neuron_astrocyte_model_wave_sol_4D_simple_eps1}
correspond to a singular perturbation problem in parameter $\varepsilon$. First, we are going to construct an heteroclinic orbit in the singular limit, that is when $\varepsilon \rightarrow 0$ by splitting the dynamics into the slow and fast dynamics. 

\subsection{Dynamics of the slow subsystem}

    Taking $\varepsilon = 0$ in system \eqref{neuron_astrocyte_model_wave_sol_4D_simple_eps}, we obtain the so-called \textit{slow subsystem} - a system of differential-algebraic equations - given by
    \begin{subequations} \label{slow_subsystem}
        \begin{align}
            0 &= \frac{1}{c} f(x,z)\, , \label{slow_subsystem1} \\[0.2em]
            0 &= \frac{1}{c} g(y,z) \, , \label{slow_subsystem2} \\[0.2em]
            \frac{d z}{d \xi} &= w \, , \label{slow_subsystem3} \\[0.2em]
            \frac{d w}{d \xi} &=  c w - h(x,y,z)\, . \label{slow_subsystem4}
        \end{align}
    \end{subequations}
    
    The first two algebraic equations \eqref{slow_subsystem1}-\eqref{slow_subsystem2} describe the commonly named \textit{critical manifold} 
    $\mathcal{M}_0$, a 2-dimensional manifold defined as  
    \begin{equation}\label{eq:crit_man}
    \begin{array}{rl}
    \mathcal{M}_0 &:= \{ (x,y,z,w) \in \mathbb{R}^4 \, | \, f(x,z)=0, \, g(y,z)=0\} \\[0.2em]
     &:= \widetilde{\M}_0 \times \mathbb{R},
     \end{array}
    \end{equation}
where 
\begin{equation}\label{eq:Mtilde0}
\widetilde{\M}_0:=\{ (x,y,z) \in \mathbb{R}^3 \, | \, f(x,z)=0, \, g(y,z)=0\}.
\end{equation}

\begin{figure}[htbp!]
        \centering
        \includegraphics[width=\textwidth]{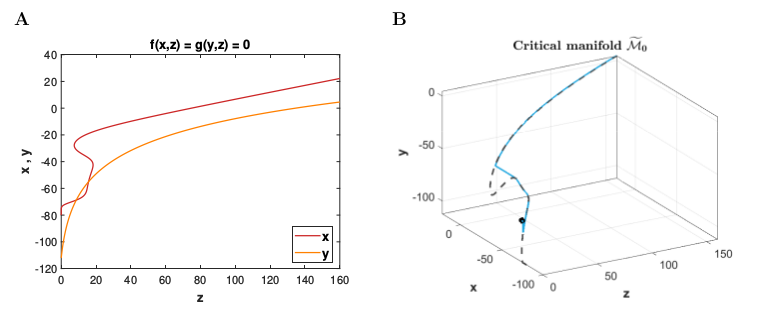}
        \caption{Dissection of the critical manifold defined by \eqref{slow_subsystem1}-\eqref{slow_subsystem2}. (A) Zero curves of $f(x,z)$ (red) and $g(y,z)$ (orange) on the $(x,z)$ and $(y,z)$ planes, respectively. (B) The manifold $\widetilde \M_{0}$ (dashed grey curve) defined in \eqref{eq:Mtilde0}, which is the projection of the 2D critical manifold $\M_{0}$ onto the $(z,x,y)$ space, and a representative solution of the simulation of system \eqref{current_pde_model} displayed in Figure \ref{simulations_rlee_reduced_gh_model} (cyan curve starting at the black-filled dot). 
        }
        \label{critical_manifold_dissection}
    \end{figure}

    Figure \ref{critical_manifold_dissection}A shows the curves $f(x,z) = 0$ and $g(y,z) = 0$ on the $(x,z)$ and $(y,z)$-plane, respectively. 
    Notice that for a finite range of values of $z$, the equation $f(x,z)=0$ has three branches $x=X^{l,m,r}(z)$ where, when they are defined, 
    \[X^l(z) \leq X^{m}(z) \leq X^{r}(z).\]
    If we denote by $(x^{L},z^{L})$ and $(x^{R},z^{R})$, the two fold points (i.e. solutions of the system $f(x,z)=0$, $f_x(x,z)=0$), such that $z^L<z^R$, then  
$X^l$ is defined for $z<z^R$, $X^m$ for $z^L<z<z^R$ and $X^r$ for $z>z^L$. 
    Notice though that the curve $g(y,z)=0$ satisfies $g_y(y,z) <0$ for all $z$ and therefore defines $y=Y(z)$ globally (see Figure \ref{critical_manifold_dissection}A). 
    Thus, the manifold $\mathcal{M}_0$ has three branches (see Figure \ref{critical_manifold_dissection}B), namely 
    \begin{equation} \label{eq:branches} 
    \begin{array}{rl}
    \mathcal{M}_0^l & =\{(x,y,z,w) \in \mathbb{R}^4 \, | \, z < z^R, x=X^{l}(z), y=Y(z)\}, \\[0.2em]
    \mathcal{M}_0^m & =\{(x,y,z,w) \in \mathbb{R}^4 \, | \, z^L < z < z^R, x=X^{m}(z), y=Y(z)\}, \\[0.2em]
    \mathcal{M}_0^r & =\{(x,y,z,w) \in \mathbb{R}^4 \, | \, z > z^L, x=X^{r}(z), y=Y(z)\}, \\
    \end{array}
    \end{equation}
as well as two fold lines $\mathcal{F}^{L,R}=\{(x^{L,R},Y(z^{L,R}),z^{L,R})\} \times \mathbb{R}$.

    In Figure \ref{critical_manifold_dissection}B we show the manifold $\widetilde{\M}_0$ defined in \eqref{eq:Mtilde0}, which is the projection  of the critical manifold $\mathcal{M}_0$ defined in \eqref{eq:crit_man} onto the $(x,y,z)$ space (dashed grey curve) and a solution of the simulation of system \eqref{current_pde_model} portrayed in Figure \ref{simulations_rlee_reduced_gh_model} for a fixed value of the spatial variable (cyan curve, starting at the black dot). This illustrates that the solution stays close to the critical manifold $\mathcal{M}_0$ whenever possible, with a rapid transition between the lower branch $\M_0^l$ and the upper branch  $\mathcal{M}^r_0$, that occurs at the fold line $\mathcal{F}^R$. So, it is reasonable to start the study by computing an approximation of the heteroclinic orbit in the singular limit that consists of pieces of solutions of the slow and fast subsystems. This will be done in Section \ref{sec0:sing_het}.

    The next step is to study the dynamics of the slow subsystem \eqref{slow_subsystem}, which describes the motion on $\M_0$. 
    The dynamics on the branch $\mathcal M^{\ast}_0$ for $\ast=l,m,r$, of the manifold $\mathcal M_0$ is given by system
    \begin{equation}\label{eq:slow_v2}
 \left \{ \begin{array}{rl}
 \dfrac{dz}{d \xi} &= w \, , \\[0.3em]
 \dfrac{d w }{d \xi} &= c w- h(X^{\ast}(z),Y(z),z)  \, ,
\end{array} \right .
\end{equation}
    where $X^{\ast}(z)$ and $Y(z)$ are the functions defining the manifolds $\M_0^{\ast}$ introduced in \eqref{eq:branches}.
    
\begin{figure}[htbp!]
        \centering
        \includegraphics[width=0.65\textwidth]{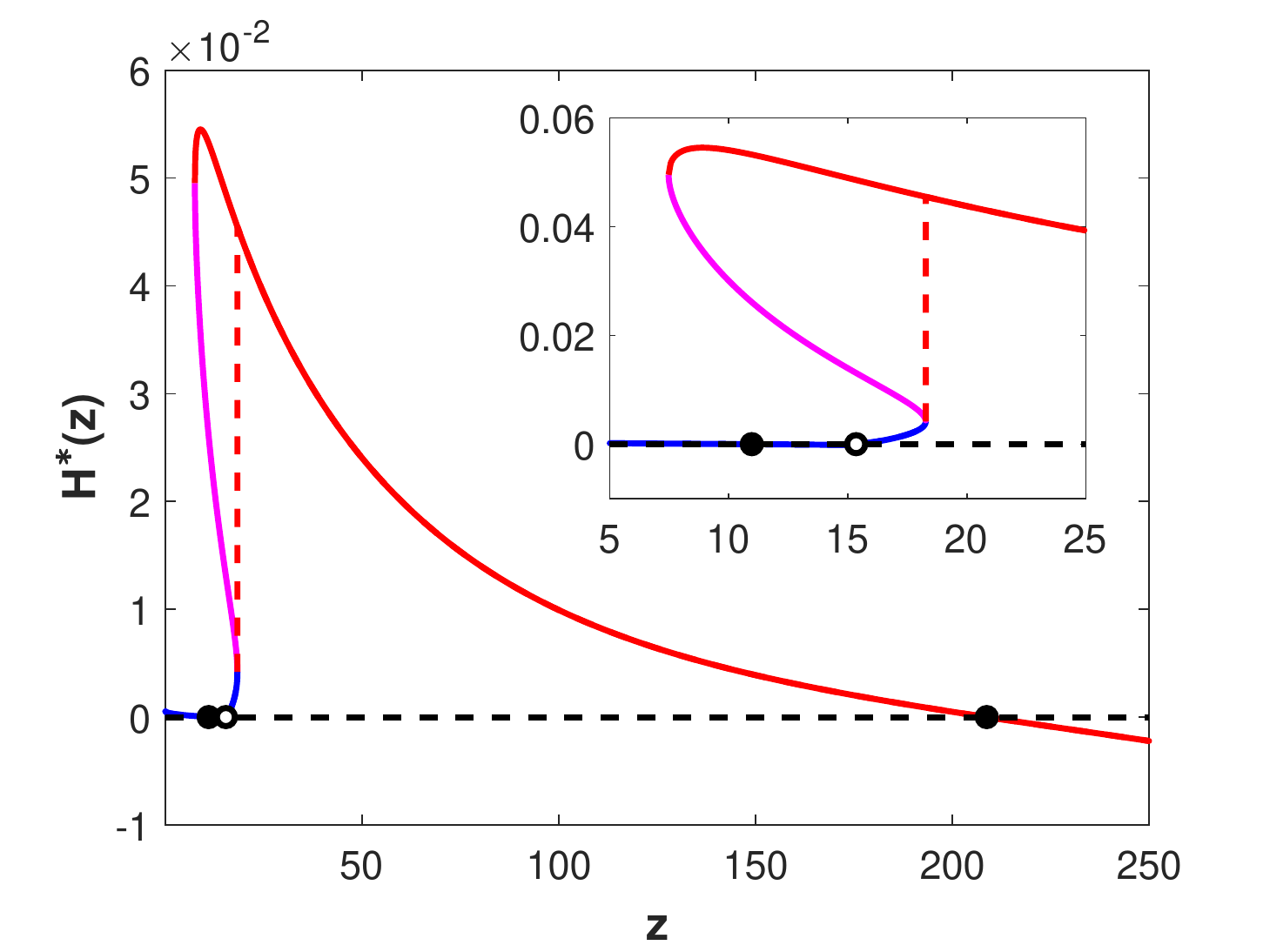}
        \caption{Graph of functions $H^{\ast}(z) := h\big(X^{\ast}(z),Y(z),z\big)$ for $\ast = l$ (blue), $m$ (purple) and $r$ (red).
        The filled (resp. empty) dots indicate the zeros of the functions $H^{\ast}$ corresponding to stable (resp. unstable) points. Inset shows a zoom of the region close to the zeros $z_1^l$ and $z_2^l$ (lying very close to each other). To compute them we have applied a Newton method to solve the system $(f,g,h)(x,y,z) = 0$. The expressions for the first partial derivatives are given in \eqref{derivatives_functions_f_g} and \eqref{derivatives_function_h}. }
        \label{equilibrium_points_slow_system}
    \end{figure}

    We begin by looking at the equilibrium points of system \eqref{eq:slow_v2}.    
    Therefore, we seek for zeros of 
    \begin{equation}\label{eq:H*}
        H^{\ast}(z):=h\big(X^{\ast}(z),Y(z),z\big) \, ,
    \end{equation} 
    for $\ast=l,m,r$. Notice that each function $H^{\ast}$ is defined for a different range of $z$. In Figure \ref{equilibrium_points_slow_system} we plot the functions $H^{\ast}(z)$ (distinguished by the colour), each one defined in its respective domain. It turns out that $H^l$ has two zeros, $z_1^l$ and $z_2^l$, $H^r$ only one, $z^r$, while $H^m$ none.  

    Notice that the equilibrium points for the slow subsystem coincide with the equilibrium points of the full system \eqref{neuron_astrocyte_model_wave_sol_4D_simple} even for $\varepsilon \neq 0$.
    The equilibrium points are given in Table \ref{tab:eq_points}. Two of the equilibrium points, $p_{l_1}=(X^l(z_1^l), Y(z_1^l),z_1^l,0)$ and $p_{l_2}=(X^l(z_2^l), Y(z_2^l),z_2^l,0)$, lie on the lower branch $\mathcal{M}_0^l$ of the manifold and the third one, $p_r = (X^r(z^r), Y(z^r),z^r,0)$, on the upper branch $\mathcal{M}_0^r$.

    \begin{table}[htbp!]
    \begin{center}
        \begin{tabular}{c|c|c|c}
             & $p_{l_1}$ & $p_{l_2}$  & $p_r$  \\
             \hline
            $x$ & -67.353771012452825 & -57.045796241401931 & 35.198894535488229 \\
            $y$ & -63.416145863486385 & -55.561014704831557 & 11.631018842324311 \\
            $z$ & 10.966529992012319 & 15.351285610517010 & 208.7014642903386
        \end{tabular}
    \end{center}
    \caption{Equilibrium points of the full system \eqref{neuron_astrocyte_model_wave_sol_4D_simple}. The component $w$ is 0.}
    \label{tab:eq_points}
    \end{table}

\begin{figure}[H]
        \centering
        \includegraphics[width=0.65\textwidth]{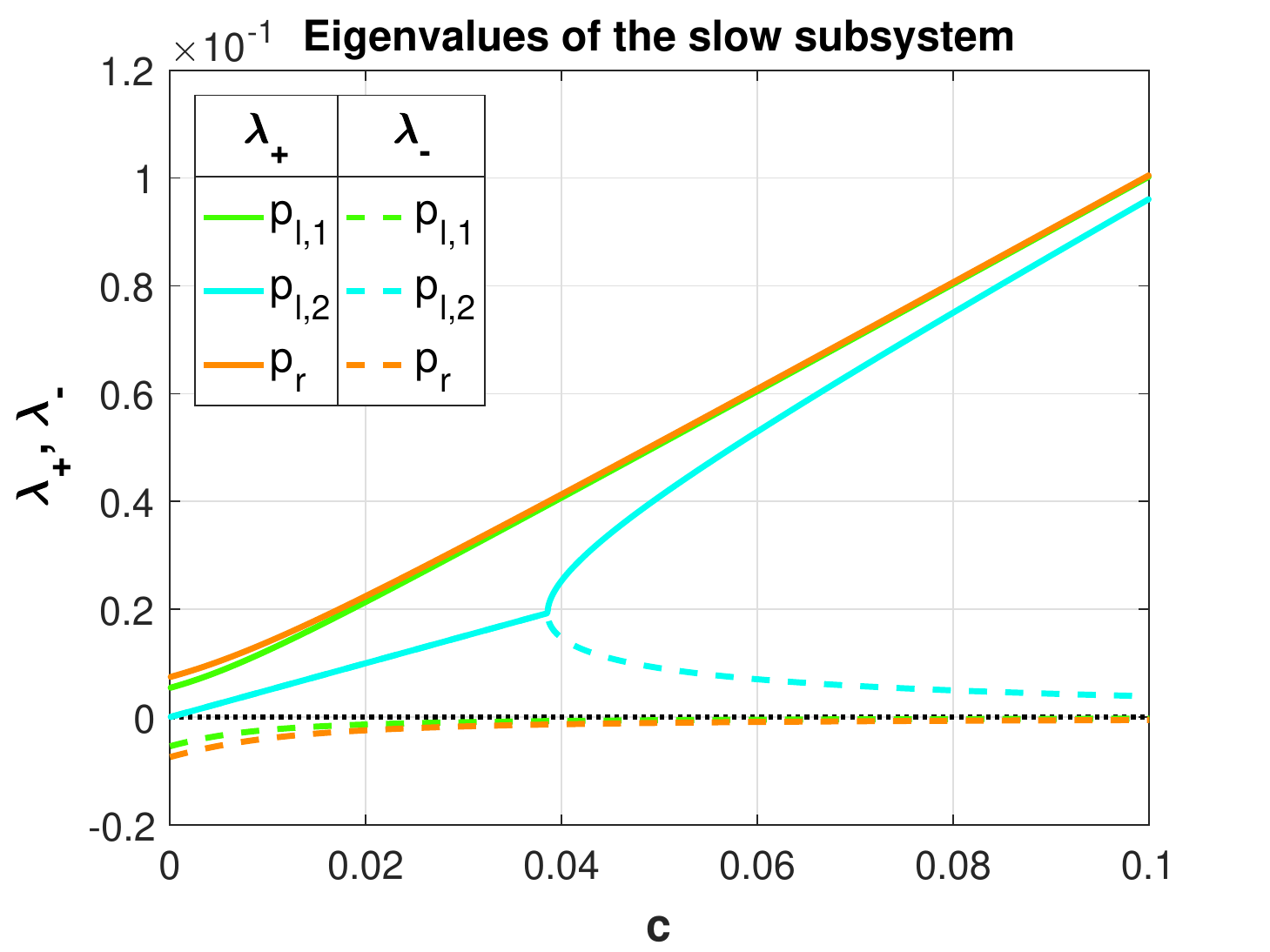}
        \caption{Eigenvalues of the slow subsystem \eqref{slow_subsystem} associated to each of its equilibrium points as $c$ varies from 0 to 0.1. Each color corresponds to an equilibrium point and the eigenvalues are differentiated by the line-style: solid lines for $\lambda_+$ and dashed ones for $\lambda_-$. Eigenvalues for $p_{l_2}$ are complex conjugates for small values of $c$ and we display only the real part of them. For the range of $c$ shown, points $p_{l_1}$ and $p_r$ are of saddle-type while $p_{l_2}$ is unstable. 
        }
        \label{eigenvalues_slow_subsystem_D1}
    \end{figure}

    If we denote by $F^{\ast}_{ss}(z,w)$, for $\ast=l,m,r$, the vector fields of system \eqref{eq:slow_v2}, the stability of the equilibrium points $p_{\ast}$ is determined by the eigenvalues $\lambda_{\pm}^{\ast}$ of the matrix 
    \begin{equation*}
        DF^{\ast}_{ss}(z^{\ast},0) = 
        \begin{pmatrix}
            0 & 1 \\
            - \frac{d H^{\ast}}{d z} & c
        \end{pmatrix} \, .
    \end{equation*}
    Observe that $\det DF^{\ast}_{ss}=\frac{d H^{\ast}}{d z}$ and $\textrm{Tr} \, DF^{\ast}_{ss}=c>0$. The expressions for the derivatives of $H^\ast$ are given in Appendix~\ref{ap:derH}. Function $H^{l}$ has negative slope at point $z_1^l$ and so does $H^r$ at $z^r$, hence they are saddle points. The derivative of $H^l$ at $z_2^l$ is positive, and thereby assuring its unstable character. To illustrate this, in Figure \ref{eigenvalues_slow_subsystem_D1} we plot both eigenvalues $\lambda_\pm$ for each equilibrium point for a range of positive values $c$.

    \subsection{Dynamics of the fast subsystem}

    Setting now $\varepsilon=0$ in system \eqref{neuron_astrocyte_model_wave_sol_4D_simple_eps1}, we obtain the so-called \textit{fast subsystem},
    \begin{subequations} \label{fast_subsystem}
        \begin{align}
            \frac{d x}{d \eta} &= \frac{1}{c} f(x,z) \label{fast_subsystem1} \, , \\[0.2em]
            \frac{d y}{d \eta} &= \frac{1}{c} g(y,z) \label{fast_subsystem2} \, , \\[0.2em]
            \frac{d z}{d \eta} &= 0 \label{fast_subsystem3} \, , \\[0.2em]
            \frac{d w}{d \eta} &= 0 \label{fast_subsystem4} \, .
        \end{align}
    \end{subequations}
    
\begin{figure}[htbp!]
        \centering
        \includegraphics[width=0.65\textwidth]{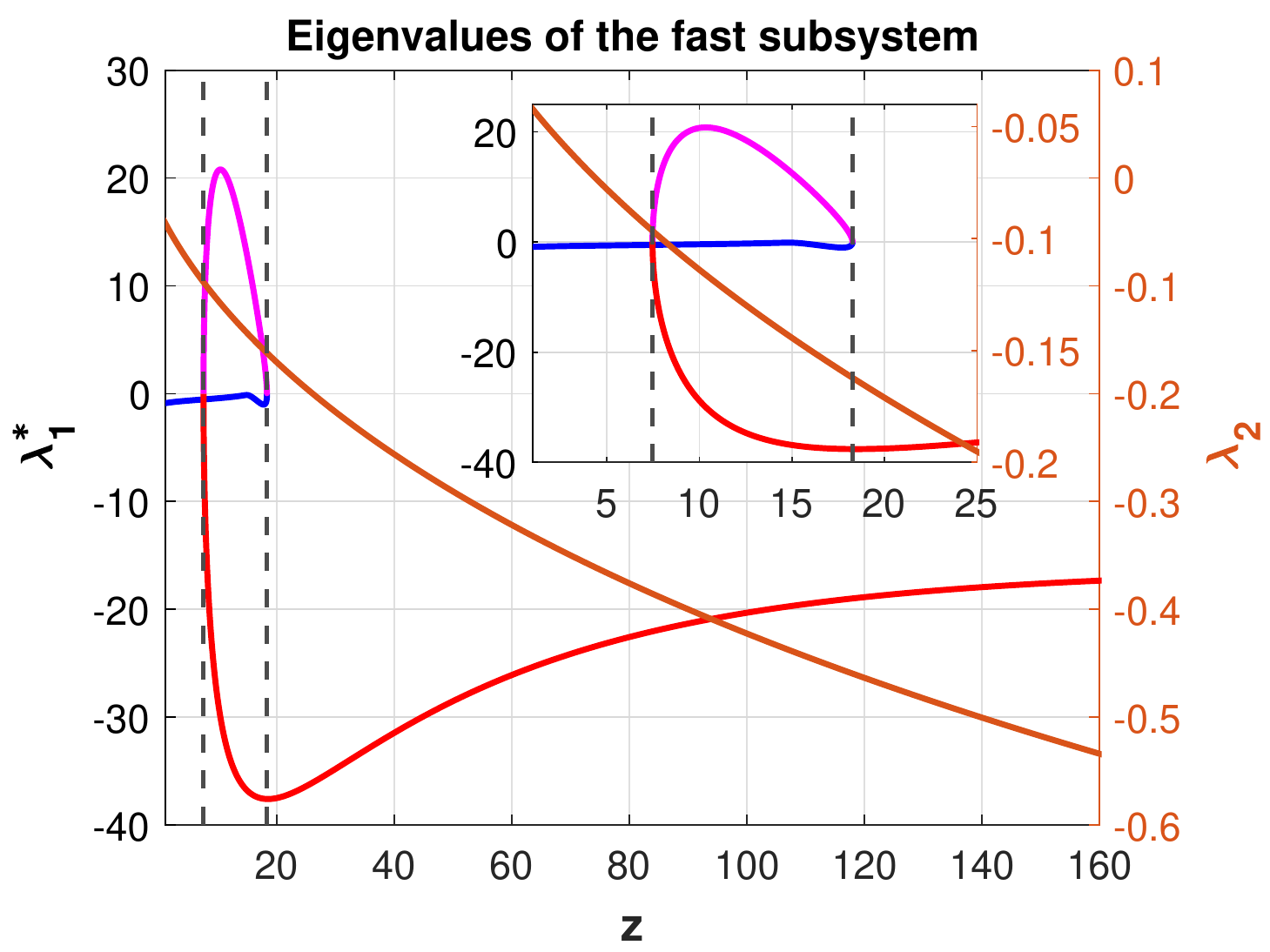}
        \caption{Eigenvalues for the equilibrium points of the fast subsystem \eqref{fast_subsystem1}-\eqref{fast_subsystem2} in terms of (parameter) $z$ for a value $c \ll 1$. On the left $y$-axis we depict the eigenvalue $\lambda_1^{\ast} (z)= \frac{1}{c} \, f_x(X^{*}(z),z)$ for each branch $X^{\ast}(z)$, for $\ast=l,m,r$ \,: $\lambda_1^l$ (blue), $\lambda_1^m$ (purple) and $\lambda_1^{r}$ (red), while on the right axis, we show in orange the second eigenvalue, $\lambda_2 (z)= \frac{1}{c} \, g_y(Y(z),z)$. Vertical dashed lines indicate the $z$-coordinate of the folds of the critical manifold (see Figure \ref{critical_manifold_dissection}). Inset shows a close-up of the region comprised by the folds, where the three solutions $X^{*}(z)$, for $\ast=l,m,r$ coexist. }\label{eigenvalues_fast_subsystem}
    \end{figure}

    We will use the fast subsystem to study the stability of the critical manifold $\mathcal{M}_0$ \eqref{eq:crit_man} in its normal directions. 
    Observe that, viewed in the fast subsystem, the critical manifold is actually a manifold filled with equilibrium points. The eigenvalues for the equilibrium points $(X^{l,m,r}(z),Y(z),z,w) \in \mathbb{R}^4$ associated to the normal directions of the critical manifold are $\lambda_1^{\ast}= \frac{1}{c} \, f_x(X^{\ast}(z),z)$ with $\ast=l,m,r$ and $\lambda_2=\frac{1}{c} \, g_y(Y(z),z)<0$. 
    In Figure \ref{eigenvalues_fast_subsystem} we plot these eigenvalues as functions of $z$.

    Observe that $\lambda_2<0$ for all values of $z$. For points on $\M_0^l$, the eigenvalue $\lambda_1^{l}<0$ and therefore, both eigenvalues are negative. The same happens for points on $\M_0^r$. However, for points on $\M_0^m$, $\lambda_1^m>0$. From this analysis we conclude that the lower and upper branches $\mathcal{M}_0^l$ and $\mathcal{M}_0^r$ are normally hyperbolic attracting manifolds, while the middle branch $\mathcal{M}_0^m$ is a normally hyperbolic saddle-type one.
    
    On a final note, the stability of the fast subsystem's equilibrium points (and hence that of the critical manifold's branches) is reversed when the wave speed $c$ becomes negative: the lower and upper branches are then unstable and the middle one is of saddle type.

    \subsection{Singular heteroclinic connection concatenating slow and fast subsystems}\label{sec0:sing_het}

    In this section we look for an heteroclinic connection between points $p_{l_1}$ and $p_r$ in the singular limit combining solutions from the slow  and fast subsytems \eqref{slow_subsystem} and \eqref{fast_subsystem}.
    Since we are looking for an orbit that tracks the critical manifold branch $\M_0^l$ until it reaches the fold $\mathcal{F}^R$, where it rapidly switches to track the upper branch $\M_0^r$, we build a piecewise system defined on the $(z,w)$ plane of the form
\begin{equation}\label{eq:piecewise}
 \begin{array}{rl}
 \dfrac{d z}{d \xi} &= w \, , \\[0.2em]
 \dfrac{d w}{d \xi} &= c w- H(z) \, ,
\end{array}
\end{equation}
    where $H(z)$ is a piecewise smooth function defined as
    \[
    H(z)= \left \{
    \begin{array}{ll}
    H^l(z) = h(X^{l}(z), Y(z),z) &\quad \text{if} \quad z \leq z^R \, , \\ 
    H^r(z) = h(X^{r}(z), Y(z),z) & \quad \text{if} \quad z \geq z^R \, ,
    \end{array}
    \right .
    \]
    where $z^R=18.276$ is the $z$-coordinate of the right fold $\mathcal{F}^R$.

    We look for an heteroclinic orbit  between the two saddle points, $\widehat p_{l_1}=(z_1^l,0)$ and $\widehat p_r=(z^r,0)$, in the two-dimensional (piecewise) system \eqref{eq:piecewise}.
    To that end, we must compute the unstable and stable invariant manifolds of $\widehat p_{l_1}$ and $\widehat p_r$, respectively, which are 1D curves, and extend them until reaching the section 
    \begin{equation}\label{eq:sigma}
    \Sigma:=\{(z,w) \in \mathbb{R}^2 \, | \, z = z^R\},    
    \end{equation} 
    where the definition of function $H(z)$ changes from the lower branch $\M_0^l$ to the upper branch $\M_0^r$ of the critical manifold. Generically, these pieces of curves will not meet $\Sigma$ at the same point. For this to occur, we treat $c$ as a parameter and seek for a value $c_0$ such that the distance between the intersection points of both invariant manifolds on $\Sigma$ is zero. Notice that the resulting heteroclinic orbit will be piecewise-defined.

\begin{figure}[htbp!]
        \centering
        \includegraphics[width=0.65\textwidth]{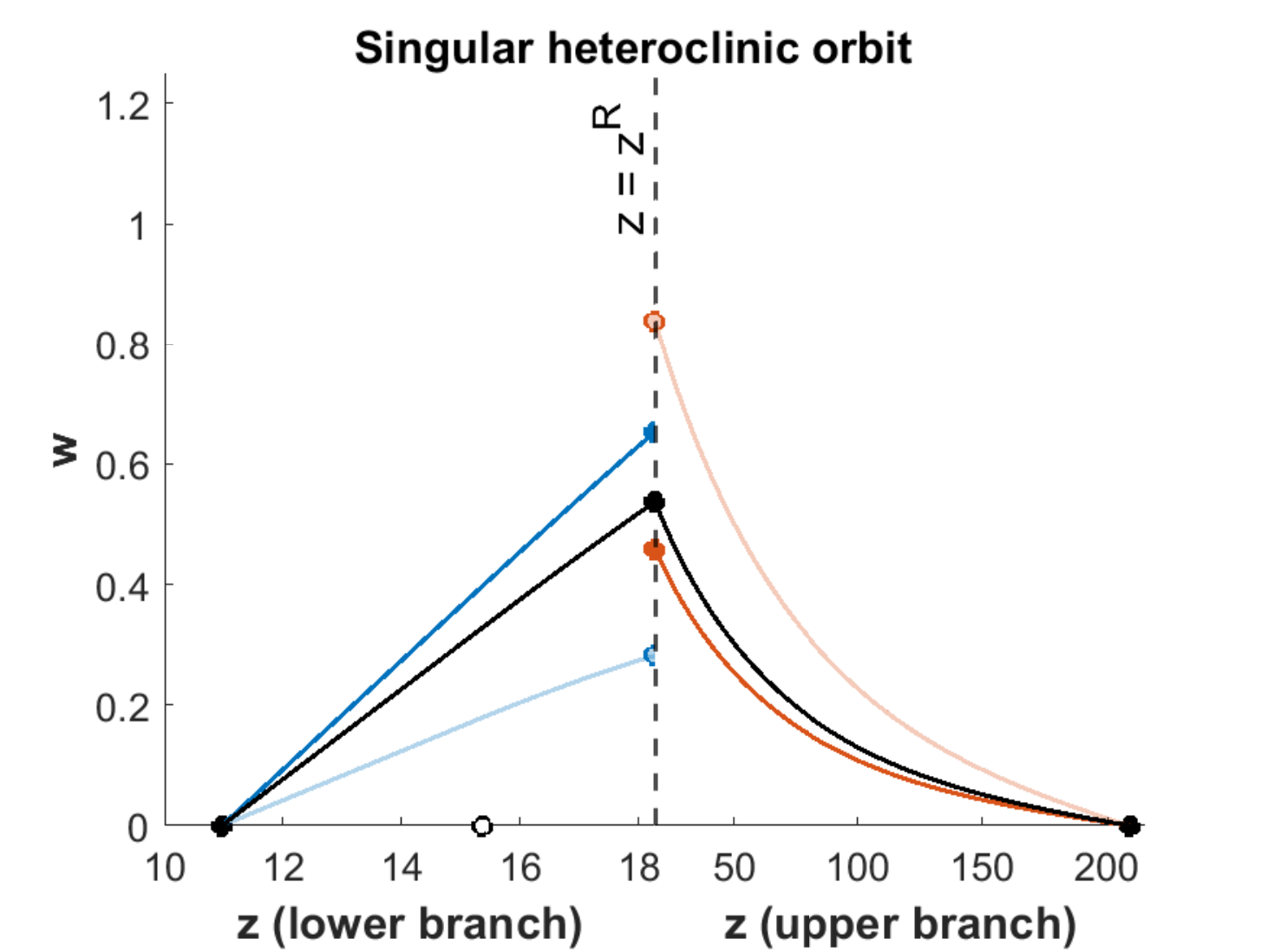}
        \caption{
        Computation of the heteroclinic orbit in the (singular) slow subsystem  \eqref{eq:piecewise}. Blue curves depict the unstable invariant manifolds of $p_{l_1}$ for two different values of $c$ up to the Poincaré section $\Sigma$, located at fold $z^R = 18.276$ (dashed-grey line). Orange curves represent the stable invariant manifolds of $p_r$ for the same two values of $c$. The curves are distinguished by the degree of transparency; the less transparent line corresponds to $c = 0.04$ while the opaque one to $c = 0.09$. The singular heteroclinic orbit (black curve) happens for $c_0 \approx 0.07426$. Notice the different ranges of the left and right intervals. This has been done for visualization purposes, so that the first part of the heteroclinic orbit can be fully appreciated.
        }
        \label{singular_heteroclinic_orbit}
    \end{figure}

    In Figure \ref{singular_heteroclinic_orbit} we plot, for different values of $c$, the unstable and stable invariant manifolds of $\widehat p_{l_1 }$ and $\widehat p_r$, respectively (blue and orange curves). The invariant manifolds have been computed using the first order approximation; we take as initial conditions $\widehat p_{\ast} + s \, v_{\ast}$ for $\ast=l_1,r$, with $v_{\ast}$ the eigenvectors of the linearization of the vector field \eqref{eq:piecewise} at $\widehat p_{\ast}$ and $s$ a small parameter, and integrate forward or backward in time accordingly. The (piecewise) heteroclinic orbit, plotted as solid black curve, is the solution of solving $d(c) := w_u - w_s = 0$, where $w_u$ and $w_s$ denote the $w$-coordinates of the intersection points of the unstable and stable manifolds with the section $\Sigma$, respectively. The singular heteroclinic orbit is found for $c_0 \approx 0.07426$.

    We recognize \eqref{eq:piecewise} as the bistable reaction-diffusion equation studied in \cite{Lee17}. Adapting classical methods described in \cite{KeenerS_book, ET_book}, one can show that the function $H(z)$ satisfies the hypothesis that guarantee the existence of heteroclinic orbits for system \eqref{eq:piecewise}.

    \section{Computing the heteroclinc orbit for system \eqref{neuron_astrocyte_model_wave_sol_4D_simple_eps} with $\varepsilon \neq 0$ by means of the parameterization method}\label{sec:full_system}

    In this section we compute numerically the heteroclinic orbit of the full system \eqref{neuron_astrocyte_model_wave_sol_4D_simple_eps} for $\varepsilon \neq 0$. 
    Recall that for $\varepsilon=0$ we had the critical normally hyperbolic attracting manifolds $\M_0^r$ and $\M_0^l$.
Fenichel's theory \cite{Fenichel71, Fenichel73} states that these critical manifolds persist, under small perturbations of $\varepsilon$,  to the full system \eqref{neuron_astrocyte_model_wave_sol_4D_simple_eps} as $\M_{\varepsilon}^l$ and $\M_{\varepsilon}^r$, while preserving in turn their hyperbolicity properties. Moreover, the new invariant manifolds $\M_{\varepsilon}^l$ and $\M_{\varepsilon}^r$ contain the critical points $p_{l_1}$ and $p_r$, respectively.  
    
    Recall that manifolds $\M_{\varepsilon}^l$ and $\M_{\varepsilon}^r$ attract strongly all nearby orbits because of their normal directions being fast and contractive.
    Therefore, both equilibrium points $p_{l_1}$ and $p_r$ maintain their saddle-type nature but with one expanding and three contractive directions.   
    This ensures that the heteroclinic orbit lives within the 1D unstable manifold of the point $p_{l_1}$, named $\Gamma^u$, until it arrives to the right fold, where it jumps to the upper branch and asymptotically approaches the manifold $\M_{\varepsilon}^r$ and follows the 1D submanifold $\Gamma^s \subset W^s(p_r)$, tangent to the slowest contracting direction at $p_r$ (which is exponentially close to the unique stable direction of $p_r$ when restricted to the critical manifold $\M^r_{\varepsilon}$). 
    
    Typically, detecting numerically an heteroclinic (or homoclinic) orbit involves following the unstable manifold of one equilibrium and the stable manifold of the other and determine whether they intersect on a convenient Poincaré section. Following manifolds of any dimension is however far from being a simple task. Even though in our present study the stable manifold of $p_r$ is 3-dimensional, there is no need to follow the whole invariant object because we know the heteroclinic orbit will be through the 1D submanifold $\Gamma^s$. 
    
    Nevertheless, finding numerically an heteroclinic connection in the 4D system \eqref{neuron_astrocyte_model_wave_sol_4D_simple_eps} proves to be a difficult task, due primarily to the presence of the two fast variables, $x$ and $y$, and the saddle character of the upper equilibrium point $p_r$. This prevents the usage of the same strategy used for the slow subsystem. On one side, following the 1D unstable manifold of $p_{l_1}$ $\Gamma^u$ works perfectly well until soon after the jump-up to the upper branch $\M^r_{\varepsilon}$, where the orbit escapes along the unstable direction of the point $p_r$. On the other hand, integrating backwards  $\Gamma^s$ blows up to infinity due to the presence of the two fast variables. Indeed, the fast variables compress the dynamics of the 4D space in the slow manifold very quickly, nonetheless backward integration reverses its stability, which turns $x$ and $y$ into expanding directions. Such fast expanding directions cause that orbits close to the 2D manifold $\M_{\varepsilon}^r$ explode to infinity at the slightest numerical error.

    For these reasons, following $\Gamma^s$ requires a more sophisticated method than simple backward integration from a linear approximation of the manifold. We will use the parameterization method \cite{Haro_book} in order to obtain a more accurate local approximation of $\Gamma^s$ which will be backward-integrated using an stiff solver (ode15s in Matlab).

    Next, we provide details of the parameterization method applied to our problem.
    We seek for a parametrization of the manifold $\Gamma^s$ as 
    \begin{equation}
    \begin{array}{rccl}
    W : &\mathbb{R} &\rightarrow &\mathbb{R}^4  \\
    & s & \rightarrow &W(s). \\
    \end{array}
    \end{equation}
    Simultaneously, we will also determine the internal dynamics of the manifold, governed by the ODE $\dot{s} = f(s)$ with $f : \mathbb{R} \rightarrow \mathbb{R}$. Imposing $W(s)$ to be a solution of the system $\dot{X} = F(X)$ (system \eqref{neuron_astrocyte_model_wave_sol_4D_simple_eps} in our case), we obtain the invariance equation
    \begin{equation} \label{invariance_equation_parameterization}
        F(W(s)) = DW(s) f(s) \, .
    \end{equation}
    Assume $W(s)$ and $f(s)$ can be both expressed as power series in $s$,
    \begin{equation} \label{coefficients_parameterization}
        W(s) = \sum_{k=0}^\infty W_k s^k \enspace , \quad f(s) = \sum_{k=0}^\infty f_k s^k \, ,
    \end{equation}
    where coefficients $W_k \in \mathbb{R}^4$ and $f_k \in \mathbb{R}$.   
    The first coefficient of the parametrization, $W_0$, is actually the equilibrium point (in our case $p_r$). 
    After substituting \eqref{coefficients_parameterization}, the invariance equation \eqref{invariance_equation_parameterization} is left as
    \begin{equation*}
        F\Bigg(\sum_{k = 0}^\infty W_k s^k\Bigg)  = \Bigg(\sum_{k = 1}^\infty W_k k s^{k-1}\Bigg) \Bigg(\sum_{k = 0}^\infty f_k s^k\Bigg) \, ,
    \end{equation*}
    to be solved order-by-order to determine the unknown coefficients $W_k$ and $f_k$. One easily obtains that $f_0=0$ so that equation of order 0 holds.
    For higher orders we develop the left-hand side of the above equation by Taylor expanding the vector field $F$ about the equilibrium point
    \begin{equation} \label{taylor_series_F_invariance_equation}
        \underbrace{F(p_r) + DF(p_r) \Bigg( \sum_{k \geq 1} W_k s^k \Bigg) + \frac{1}{2!} D^2 F(p_r) \Bigg( \sum_{k \geq 1} W_k s^k \Bigg)^2 + \ldots}_{\sum_{k \geq 0} F_k s^k} = \Bigg(\sum_{k = 1}^\infty W_k k s^{k-1}\Bigg) \Bigg(\sum_{k = 1}^\infty f_k s^k\Bigg) \, .
    \end{equation}
    The equation of order 1 leads to the eigenvalue problem
    \begin{equation}
        DF({p_r}) W_1 = f_1 W_1 \, ,
    \end{equation}
    thus begin $W_1$ the eigenvector of eigenvalue $f_1$. Since we are looking for the slow (sub)manifold (associated to the slowest eigenvalue), we choose $f_1 := \lambda_{\text{slow}}$ and $W_1$ the corresponding eigenvector, that we will take of modulus 1.

    Solving the equations of higher order, shows that there is not a unique solution. We will use this freedom to choose the dynamics $f$ to be as simple as possible, namely, as
    \begin{equation}\label{eq:f(s)}
        f(s)=\lambda_\text{slow} s,
    \end{equation}
    which corresponds to $f_k=0$ for $k \geq 2$.

    With this simplification, the equation of order 2 reads as
    \begin{equation}
        \underbrace{ \, DF(p_r) W_2 + \big[ F(W_{< 2}) \big]_{|_2} \, }_{F_2} = 2 \lambda_{\text{slow}} W_2 \, ,
    \end{equation}
    where $W_{<2}=\sum_{k < 2} W_k s^k$ and $\big[ F(W_{< 2}) \big]_{|_2}$ denotes the coefficient of order 2 of the Taylor expansion in $s$. 
    The resulting system is linear
    \begin{equation}
        (DF(p_r) - 2 \lambda_{\text{slow}} \textrm{Id}) W_2 = -\big[ F(W_{< 2}) \big]_{|_2} \, ,
    \end{equation}
    and therefore it can be solved efficiently. In general, for the coefficient of order $k \geq 2$, $W_k$, we need to solve the linear system
    \begin{equation} \label{param_method_linear_system_orderk}
        (DF(p_r) - k \lambda_{\text{slow}} \textrm{Id}) W_k = -\big[ F(W_{< k}) \big]_{|_k} \, ,
    \end{equation}
    where, again, the coefficient $F_k$ has been split into the known, $\big[ F(W_{< k}) \big]_{|_k}$, and unknown parts, $DF(p_r) W_k$.

    Equation \eqref{param_method_linear_system_orderk} can be solved provided 
    \[ \det (DF(p_r) - k \lambda_{\text{slow}} \, \textrm{Id}) \neq 0 \, ,
    \]
    or, equivalently, that $k \lambda_{\text{slow}}$ is not an eigenvalue of the differential matrix $DF(p_r)$, which is fulfilled in our case because the other eigenvalues are much larger  
    (this is equivalent to say there are no resonances among the eigenvalues of the system). 
    
    The right-hand side of equation \eqref{param_method_linear_system_orderk} can be numerically computed efficiently with accuracy as high as desired using the Automatic Differentiation algorithms \cite{Haro_book, Griewank_book}. The algorithm outputs the coefficients (up to any order) of the power series of the composition of a truncated power series with any analytic function $F$. The method consists in computing systematically compositions of power series with elementary functions, for which there are well-established recurrences (see for instance \cite{Haro_book}).
    
\begin{figure}[htbp!]
        \centering
        \includegraphics[width=0.48\textwidth]{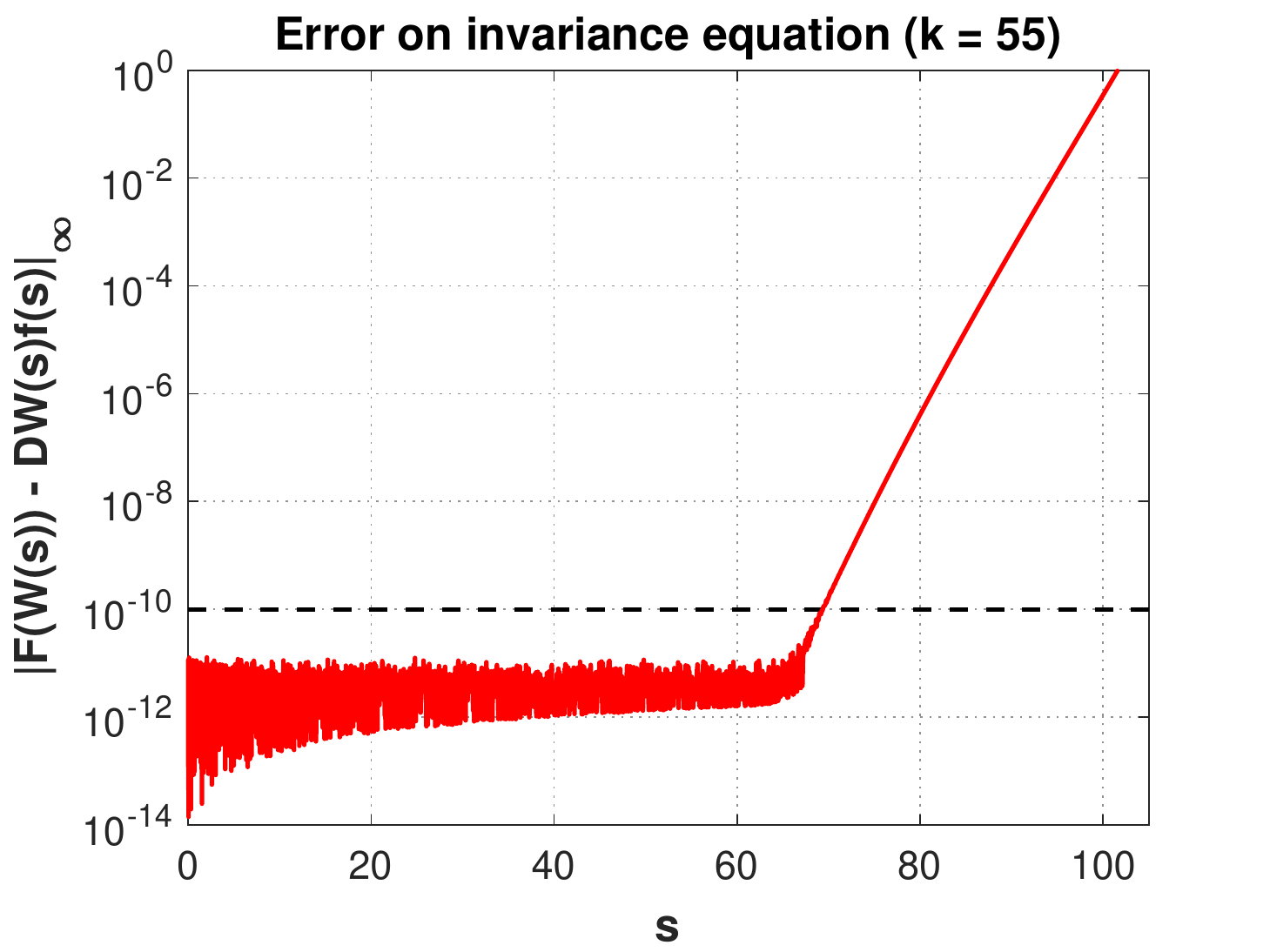}
        \caption{ Error of the invariance equation \eqref{invariance_equation_parameterization} as a function of the parameter $s$ for a parametrization $W(s)$ computed up to order $k=55$ with the internal dynamics $f(s)$ given by \eqref{eq:f(s)}. Any value $s$ whose error is smaller than the allowable threshold (i.e. below the dashed-black line) produces a potentially accurate initial condition to be backward-integrated.}
        \label{error_invariance_equation_c006}
    \end{figure}

    Once fixed the order of the expansion $W(s)$ and given its coefficients $W_k$, we take the largest value $s$  
    such that the error in the invariance equation, $||F(W(s)) - DW(s) f(s)||_{\infty}$, is less than $10^{-10}$ (see Figure \ref{error_invariance_equation_c006}) and the backward integration from $W(s)$ provides a trajectory that hits the Poincaré section $\Sigma_1:=\{z=22\}$.

    To obtain the 1D unstable manifold of the point $p_{l_1}$, $\Gamma^u$ we integrate the system forward starting from a point obtained using a first order approximation of the manifold (i.e. $p_{l_1} + \delta \, v_{l_1}$ with $v_{l_1}$ the eigenvector of the linearized system around $p_{l_1}$ associated to the positive eigenvalue and $\delta$ a small value).

\begin{figure}[htbp!]
 \centering
        \includegraphics[width=0.9\textwidth]{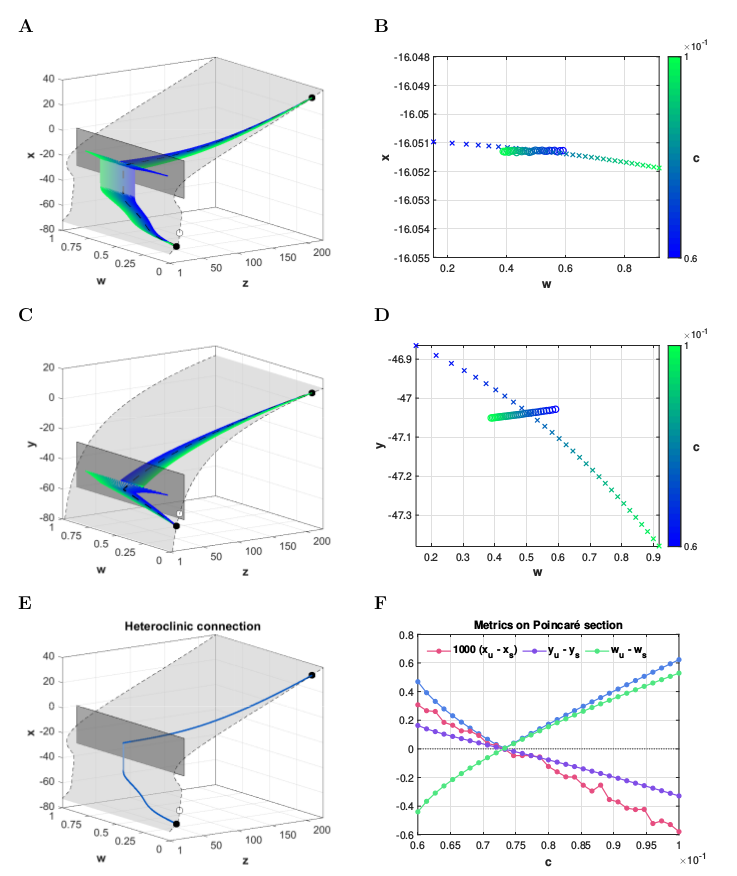}
        \caption{\footnotesize Computation of an heteroclinic orbit between the two saddle points $p_{l_1}$ and $p_r$ in the full system \eqref{neuron_astrocyte_model_wave_sol_4D_simple_eps} using the parameterization method. (A,C) 3D projections of the phase space onto the: (A) $(x,z,w)$ and (C) $(y,z,w)$ coordinated systems. In light grey we depict the critical manifold $\M_{0}$ and in dark grey we represent the Poincaré section $\Sigma_1=\{z = 22\}$. The equilibrium points are denoted by circles, of which the black-filled ones correspond to the saddle points $p_{l_1}$ and $p_r$. The coloured orbits correspond to the unstable $\Gamma^u$ and (weakest) stable $\Gamma^s$ 1D manifolds of the saddle points for different values of $c$ (colour palette). The 3D representation of the singular heteroclinic orbit, shown in Figure \ref{singular_heteroclinic_orbit}, has been represented by a black-dashed curve. (B,D) 2D projections of intersection points of the 1D manifolds $\Gamma^u$ (crosses) and $\Gamma^s$ (circles) with the Poincaré section $\Sigma_1$ for the range of values $c$ considered. Projection onto the: (B) $(w,x)$ and (D) $(w,y)$ planes. (E) 3D projection (onto the $(z,w,x)$-space) of the numerically computed heteroclinic connection for $\widehat c$. (F) Metrics computed upon the Poincaré section between the intersection points of the unstable and stable manifolds as a function of the parameter $c$: Euclidean distance (blue-dotted curve) and component-wise distances (see panel legend).}
        \label{heteroclinic_orbit_full_system}
    \end{figure}
    
    In Figure \ref{heteroclinic_orbit_full_system}A-D we illustrate the strategy used for computing an heteroclinic orbit. For a range of $c$-values including the value $c_0 = 0.07426$ (which is the value for which we found an heteroclinic orbit in the singular case), we track both 1D invariant manifolds until they reach the Poincaré section $\Sigma_1$ (see the 3D phase space projections in Figure \ref{heteroclinic_orbit_full_system}A and C). The set of intersection points of the unstable manifold $\Gamma^u$  with the Poincaré section $\Sigma_1$ forms a curve in terms of $c$ which intersects with that of the stable manifold $\Gamma^s$ (see Figure \ref{heteroclinic_orbit_full_system}B and D).  For an heteroclinic solution to exist, such an intersection should happen for the same value $c$. We therefore check that the distance $\| \cdot \|_2$ between the intersection points of both manifolds with the section $\Sigma_1$ becomes zero for some value $\widehat c$ (see blue-dotted line in Figure \ref{heteroclinic_orbit_full_system}F). In  Figure \ref{heteroclinic_orbit_full_system}F the Euclidean distance has been disaggregated into the distances (with sign) of each of its components (see pink/purple/green-dotted lines).  By applying a secant method in the $w$-component, we found an approximate value $\widehat c = 7.3135 \cdot 10^{-2}$ with an error of order $10^{-4}$.
    This corresponds to an estimated wave velocity of $vel = \widehat c \sqrt{D_k}=  7.3135 \cdot 10^{-2} \,\text{ms}^{-1/2} \cdot \sqrt{1.96 \cdot 10^{-5} \,\text{cm}^2 \, \text{s}^{-1}} = \sqrt{1.96} \cdot 0.073135 \, \text{mm/s} = 0.102389  \, \text{mm/s} =  6.1433 \, \text{mm/min}$.

\section{Computing the heteroclinic orbit using a different strategy based on Fenichel's theory}\label{sec:fenichel}

     In this section we present a different strategy to overcome the numerical difficulties to compute the 1D slow manifold $\Gamma^s$ of the critical point $p_r$. Our idea is to compute $\Gamma^s$ by integrating backwards system \eqref{neuron_astrocyte_model_wave_sol_4D_simple_eps} starting from a local approximation of it, while restricting the trajectory to be approximately on $\M^r_{\varepsilon}$, the perturbed upper branch of the manifold $\M^r_{0}$. By doing so, we assure numerical integration will not blow up to infinity either in the $x$ or $y$ directions. As discussed in the previous section, the computation of the 1D unstable manifold of $p_{l_1}$, $\Gamma^u$, does not exhibit these problems and for this reason we focus only on the point $p_r$.

    Since the upper branch of the Fenichel's manifold $\M_{\varepsilon}^r$ perturbs from the critical manifold $\M_{0}^r$, it can be approximated by considering both fast variables $x$ and $y$ as formal expansions in $\varepsilon$, 
    \begin{equation} \label{fenichel_manifold_formal_expression}
        \begin{split}
            x &= m(z,w,\varepsilon) = m_0(z,w) + \varepsilon m_1(z,w) + \varepsilon^2 m_2(z,w) + \ldots \, , \\[0.2em]
            y &= n(z,w,\varepsilon) = n_0(z,w) + \varepsilon n_1(z,w) + \varepsilon^2 n_2(z,w) + \ldots \, .
        \end{split}
    \end{equation}

    The coefficients will be determined by solving the invariance equation order by order. 
    Actually, the coefficients of order 0 are given by $m_0(z,w)=X^{r}(z)$ and $n_0(z,w)=Y(z)$ (see equation \eqref{eq:branches}).

    From imposing the manifold \eqref{fenichel_manifold_formal_expression} to be invariant by system \eqref{neuron_astrocyte_model_wave_sol_4D_simple_eps} we reach the following system of equations:
    \begin{equation}\label{invariance_equation_xy}
        \begin{split}
            \dot{x} &= \dot{m}(z,w,\varepsilon) \quad \Longrightarrow \quad 
            f\big(m(z,w,\varepsilon),z\big) = \varepsilon c w \, \frac{\partial m}{\partial z} + \varepsilon c \big(c w - h\big(m(z,w,\varepsilon), n(z,w,\varepsilon), z\big) \big) \, \frac{\partial m}{\partial w} \, , \\[0.2em]
            \dot{y} &= \dot{n}(z,w,\varepsilon) \quad 
            \Longrightarrow \quad 
            g\big(n(z,w,\varepsilon), z\big) = \varepsilon c w \, \frac{\partial n}{\partial z} + \varepsilon c \big(c w - h\big(m(z,w,\varepsilon), n(z,w,\varepsilon),  z\big) \big) \, \frac{\partial n}{\partial w} \, .
        \end{split}
    \end{equation}

    Developing functions $f$, $g$ and $h$ as Taylor series around the coefficients of order 0, $m_0$ and $n_0$, and substituting them back into system \eqref{invariance_equation_xy}, we are able to solve the resulting system of equations up to any power of $\varepsilon$. Indeed, one can obtain the following expressions (see Appendix \ref{sec:expansionmanifold}) up to order 2:
               \[ m_1(z,w) =  -\frac{c w \, f_z(X^r(z),z)}{\big(f_x(X^r(z),z)\big)^2} \, , \qquad n_1(z,w) = -\frac{c w \, g_z(Y(z),z)}{\big(g_y(Y(z),z)\big)^2} \, ,\]
    and 
    \[m_2(z,w) = -\frac{1}{2!} \frac{f_{x x}(m_0,z)}{f_x(m_0,z)} \, m_1^2 + \frac{c w}{f_x(m_0,z)} \frac{\partial m_1}{\partial z} - c \big(c w - h(m_0,n_0,z) \big) \frac{c \, f_z(m_0,z)}{\big(f_x(m_0,z)\big)^3} \, ,
    \]
    and
    \[ n_2(z,w) = -\frac{1}{2!} \frac{g_{y y}(n_0,z)}{g_y(n_0,z)} \, n_1^2 + \frac{c w}{g_y(n_0,z)} \frac{\partial n_1}{\partial z} - c \big(c w - h(m_0,n_0,z) \big) \frac{c \, g_z(n_0,z)}{\big(g_y(n_0,z)\big)^3} \, .\]

    Next, we will restrict the dynamics of the full system \eqref{neuron_astrocyte_model_wave_sol_4D_simple_eps} onto the second order approximation of the manifold $\M^r_{\varepsilon}$ by taking 
    \begin{equation}\label{eq:embed}
     \begin{array}{ll}
    x&=X^r(z)+ \varepsilon m_1 (z,w) + \varepsilon^2 m_2(z,w) \, , \\
     y&=Y(z)+ \varepsilon n_1 (z,w) + \varepsilon^2 n_2(z,w) \, . \\ 
     \end{array}   
    \end{equation}
    The resulting system is
     \begin{equation} \label{slow_subsystem_like}
        \begin{split}
            \dot{z} &= \varepsilon w \, , \\[0.2em]
            \dot{w} &= \varepsilon \big(c w - h(X^r(z)+ \varepsilon m_1 (z,w) + \varepsilon^2 m_2(z,w), Y(z)+ \varepsilon n_1 (z,w) + \varepsilon^2 n_2(z,w), z) \big) \, .
        \end{split}
    \end{equation}

    \begin{figure}[htbp!]
        \centering
        \includegraphics[width=0.9\textwidth]{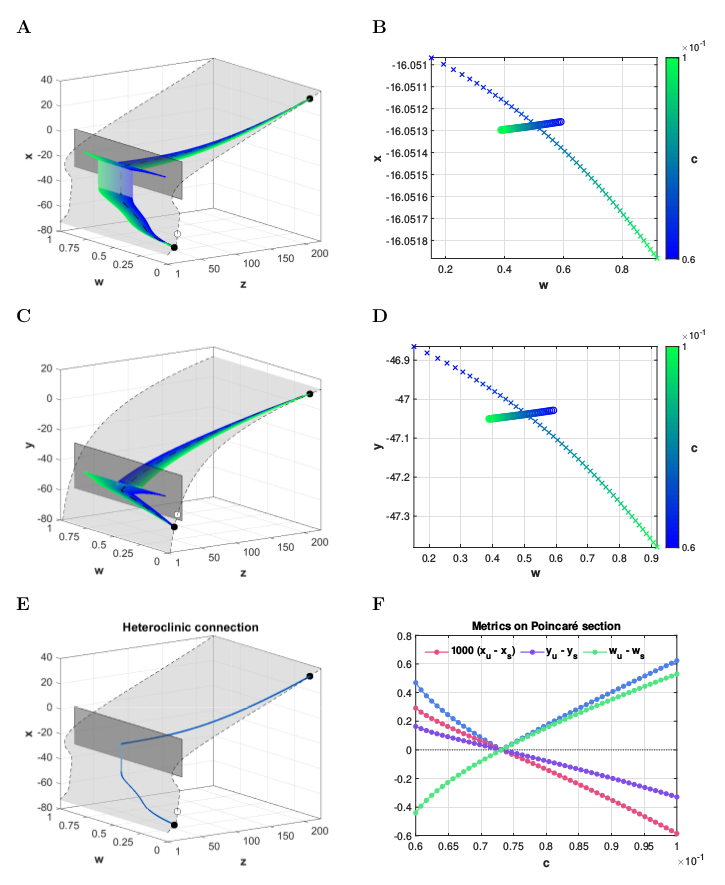}
        \caption{\footnotesize Computation of the heteroclinic orbit between the two saddle points $p_{l_1}$ and $p_r$ when full system \eqref{neuron_astrocyte_model_wave_sol_4D_simple_eps} is restricted to the Fenichel's manifold \eqref{fenichel_manifold_formal_expression} on the upper branch. (A,C) 3D projections of the phase space onto the: (A) $(z,w,x)$ and (C) $(z,w,y)$ coordinated systems. The critical manifold $\M_{0}$ and the Poincaré section $\Sigma_1=\{z = 22\}$ are depicted in light and dark gray colours, respectively. Black-filled (resp. empty) dots indicate the saddle (resp. unstable) equilibrium points of the system. The coloured curves correspond to the unstable $\Gamma^u$ and (weakest) stable $\Gamma^s$ 1D manifolds of the saddle points (up to the Poincaré section) for different values $c$ ranging from 0.06 (blue) to 0.1 (green). The black-dashed curve illustrates the 3D representation of the singular heteroclinic orbit shown in Figure \ref{singular_heteroclinic_orbit} 
        (B,D) 2D projections of the intersection points of the 1D manifolds $\Gamma^u$ (crosses) and $\Gamma^s$ (circles) with the Poincaré section $\Sigma_1$ for the range of values $c$ considered. Projection onto the: (B) $(w,x)$ and (D) $(w,y)$ planes. (E) 3D projection (onto $(z,w,x)$-space) of the numerically computed heteroclinic connection for $\widetilde c$. (F) Metrics computed upon the Poincaré section between the intersection points of the unstable, $v_u$, and the stable, $v_s$, manifolds as a function of the parameter $c$: Euclidean distance (blue-dotted curve) and component-wise distances (see panel legend).}
        \label{xy_invariant}
    \end{figure}  
    We compute the 1D stable manifold $\widetilde \Gamma^s$ of the critical point $(z^r,0)$ of system \eqref{slow_subsystem_like}. We take the local approximation provided by the linear terms and we integrate it backwards in time until the trajectory hits the Poincar\'e section $\Sigma_1=\{z=22\}$. Details are given in Appendix~\ref{ap:DF}. Thus, $\Gamma^s$ can be obtained embedding $\widetilde{\Gamma}^s$ in $\mathbb{R}^4$ using \eqref{eq:embed}. In Figure \ref{xy_invariant}A and C we show manifolds $\Gamma^s$ and $\Gamma^u$ (where the latter has been computed as described in the previous section) and their intersection along the Poincaré section $\Sigma_1$ (depicted in dark grey), for different values of the parameter $c$ (color palette).
    
    The intersection points with section $\Sigma_1$ will have $x$, $y$ and $w$ components for each value of $c$ (see Figure \ref{xy_invariant}B and D), and we look for $\widetilde c$ such that the triplet $v_u := (x_u,y_u,w_u) \in \Gamma^u$ equals $v_s := (x_s,y_s,w_s) \in \Gamma^s$. We compute the euclidean norm of the difference $v_u - v_s$ and, separately, the distances with sign for each of its components as a function of the parameter $c$ (see Figure \ref{xy_invariant}F). As we can observe, the euclidean distance vanishes within the range $(0.07, 0.075)$, and using a secant method in the $w$-component we find that it vanishes for $\widetilde c = 0.073135$ with an error of order $10^{-4}$. Notice that this value coincides with the one obtained using the parameterization method and corresponds to a wave velocity of $6.14 \, \text{mm/min}$. In Figure \ref{xy_invariant}E we display the heteroclinic connection between $p_{l_1}$ and $p_r$ resulting for $\widetilde c$.

   \section{Discussion}\label{sec:discussion}

In this paper we have presented a computational and mathematical framework to study the initiation and propagation of CSD waves in a reduction of the model introduced in \cite{Huguet16}, involving the dynamics of the extracellular potassium concentration $[K^+]_e$.  More precisely, we have studied semi-analytically the existence of traveling wave solutions in model \eqref{current_pde_model}, which correspond to heteroclinic orbits of the system \eqref{neuron_astrocyte_model_wave_sol_4D}. We have combined techniques from singular perturbation theory \cite{Jones_book, Kuehn15}, and numerically efficient methods based on the parameterization method \cite{Haro_book}, involving automatic differentiation techniques, which provide an efficient computation of the heteroclinic solution and, subsequently, the propagation velocity of the wave. The methods have proven successful to overcome the difficulties associated to the presence of different time-scales in the system.

We emphasize that the computation of the propagation velocity $vel \approx c \sqrt{D_K}$ by means of computing an heteroclinic connection is more accurate since it does not have the boundary effects that appear in the simulation of the discretized PDE model \eqref{current_pde_model} when the number of neurons is small. Indeed, the estimated velocity from the heteroclinic is faster (approx 6 mm/s) than the one obtained from numerical simulations of the PDE (approx 4 mm/s). Moreover, it is more computationally efficient since avoiding the boundary effects in the full model requires to use a large number of neurons, which in turn increases the computational cost of the simulations.

\begin{figure}[htbp!]
 \centering
        \includegraphics[width=0.9\textwidth]{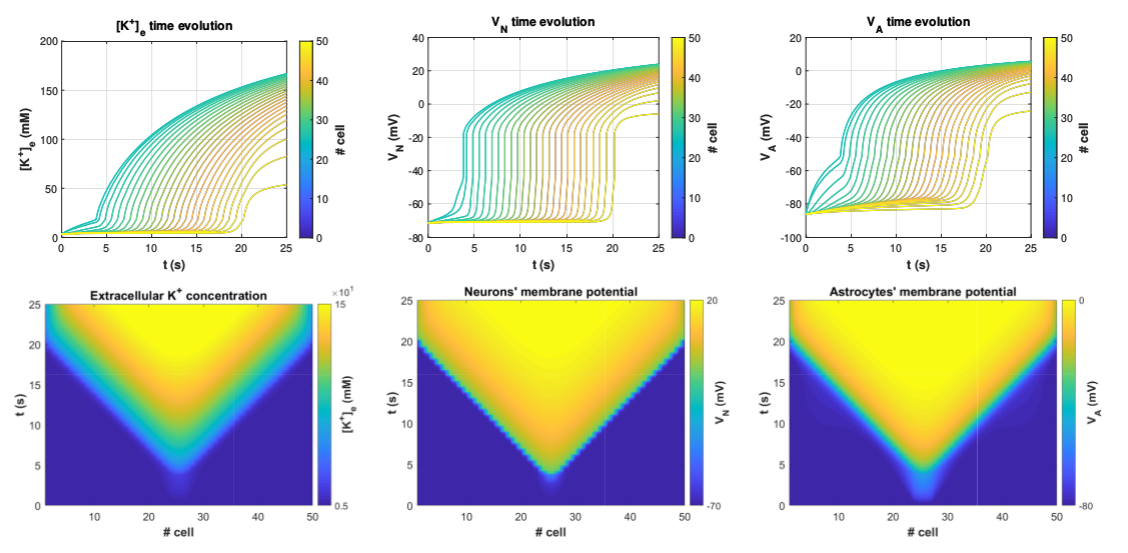}
        \caption{Simulations as in Figure \ref{simulations_rlee_reduced_gh_model} for the reduced model \eqref{current_pde_model} with $V_N$ and $V_A$ now regarded as variables with instantaneous dynamics, that is, setting $\frac{d V_N}{d t} = \frac{d V_A}{d t} = 0$ on system \eqref{current_pde_model}. Variables $V_N$ and $V_A$ are described in terms of the extracellular concentration of $K^+$ by means of the algebraic constraints given in \eqref{eq:Lee_red}. Comparing with the simulation in Figure \ref{simulations_rlee_reduced_gh_model}, the dynamics of the model \eqref{current_pde_model} is close to that of the limit case considered here.}
        \label{simulations_reduced_model_order0}
    \end{figure}

Our approach is based on the model reduction suggested in \cite{Lee17}. However, we want to highlight some relevant differences with respect to this work. In \cite{Lee17} the system is further reduced to a single PDE for $[K^+]_e$ of reaction-diffusion type, imposing that the dynamics for $V_N$ and $V_A$ in \eqref{current_pde_model} are set to be instantaneous (i.e. $dV_n/dt=dV_A/dt=0$). When imposing a traveling wave solution in this 1D reaction-diffusion PDE, it yields the dynamics described by the piecewise slow subsystem \eqref{eq:piecewise}. 
In our study, we go beyond the singular limit and compute the heteroclinic orbit and the velocity in the case $\varepsilon \neq 0$, which corresponds to the full system \eqref{current_pde_model}. We emphasize though that the singular heteroclinic orbit provides a good approximation of the dynamics since the estimated velocity is similar and the trajectory tracks closely to the singular heteroclinic orbit (see black-dashed curves in Figures \ref{heteroclinic_orbit_full_system}AC and \ref{xy_invariant}AC). 
Indeed, in Figure \ref{simulations_reduced_model_order0} we present the numerical simulations of \eqref{current_pde_model}  taking $V_N$ and $V_A$ as instantaneous variables (i.e. $\frac{d V_N}{d t} = \frac{d V_A}{d t}= 0$). As expected from our results, these simulations reveal a good agreement with those for system \eqref{current_pde_model} shown in Figure \ref{simulations_rlee_reduced_gh_model}.
  
The analysis performed herein was based on the existence of two different temporal scales in the reduced model \eqref{current_pde_model}: membrane potential variables change at a faster rate than the extracellular $K^+$ concentration. In the light of the observations in Figure \ref{functions_fgh} (see Remark \ref{rem:time_scales}), it seems reasonable to ponder another scale dynamics disparity, between the neuronal and astrocytic membrane potentials, and thus introduce a new intermediate scale $\delta$ on system \eqref{neuron_astrocyte_model_wave_sol_4D_simple_eps}. We leave this more detailed analysis for future work. 

The presence of different time-scales in a dynamical system can present a significant computational challenge. In this paper, we demonstrate the power and versatility of the parameterization method \cite{CabreFL05, Haro_book}, widely used in the field of dynamical systems \cite{HuguetL13, PerezCervera20}, to overcome some issues related to this challenge. Specifically, we applied the parameterization method to numerically compute the upper branch of the heteroclinic orbit in the system $\varepsilon \neq 0$ where conventional numerical methods failed due to the presence of a strong expansive directions.

We acknowledge that, even if our results on the existence of a traveling wave solution for the system with $\varepsilon \neq 0$ rely on accurate numerical computations, they are not a rigorous analytical proof of their existence. To do so, one needs to use techniques from singular perturbation theory \cite{Jones_book,Kuehn15} to describe the behaviour of the solutions close to the fold of the slow manifold and techniques from computer assisted proofs \cite{Capinski22} to keep track of the stable and unstable manifolds of the saddle points, until they reach suitable Poincar\'e sections, which we plan to attempt in future studies. 

   In this paper, we have assumed a reduced model \eqref{current_pde_model} without gap junctions in the astrocyte cells. However, previous work \cite{Huguet16} has shown that gap junctions play a significant role in preventing wave propagation. Therefore, future studies will involve incorporating gap junction currents and exploring how they affect the wave shape and dynamics. 

   Finally, we want to emphasize that the reduced model \eqref{current_pde_model} only includes the biophysics terms affecting the wave initiation but does not model the wave termination. Indeed, $Na^+$ concentrations are assumed constant in the reduced model, but in the original model presented in \cite{Huguet16}, the increase in $[K^+]_e$ is accompanied by a steep decline in $[Na^+]_e$ (see Figure~\ref{simulations_ghuguet_model}). Moreover, freezing the inactivation variable $h_p$ of the persistent sodium current $I_{NaP}$ prevents the termination of the wave and retains the network cells in the depolarized state.  Further analysis may focus on the impact of including the dynamics for these variables, by designing a different reduced system which, at this turn, involves only the biophysics processes involved in the wave termination. The description of the wave termination would be analogous to that conducted within this manuscript, but the heteroclinic orbit would go from the equilibrium point corresponding to the depolarized state to the one corresponding to the resting state.
   Thus, in the future we plan to concatenate two reduced systems, one that describes the wave initiation and another describing the wave termination, using a piecewise formulation \cite{RinzelK73, McKean70, KeenerS_book}.  This more accurate description of the wave will provide new insights into the dynamics of the system describing CSD, leading to a more comprehensive understanding of the mechanisms underlying the CSD phenomenon.
   
   In summary, we have presented a computational and mathematical framework that offers a more computationally efficient and accurate approach to identify the properties of CSD waves in model \eqref{gHuguet_model} compared to the crude simulation of the model. We hope that it establishes a solid foundation upon which future studies of more complex models can be built (such as extensions of the models in \cite{Lemaire23} incorporating the spatial domain). 

\newpage

\appendix

\section*{Appendixes}\label{ap:appendix}

In this section we include details of the numerical  and algebraic mathematical computations.

\section{Solving numerically the PDE systems}\label{ap:pde}
    In this section we discuss the numerical details for the simulations of the PDE systems \eqref{gHuguet_model}, \eqref{current_pde_model} and \eqref{current_pde_model} with the reductions in \eqref{eq:Lee_red} (Figures \ref{simulations_ghuguet_model}, \ref{simulations_rlee_reduced_gh_model} and \ref{simulations_reduced_model_order0}, respectively).
We use a network of N neuron-astrocyte pairs sharing the same extracellular space. Thus, the space is discretized taking a grid $x_i=i \Delta x$, for $i=1,\ldots, N$ and $\Delta x$ is the distance between neurons.
    
    The diffusion terms  are approximated by central second-order finite differences, for instance, in the $[K^+]_e$ case,
    \begin{equation}
        \frac{\partial^2 [K^+]_e}{\partial x^2} \approx \frac{[K^+]_e^{i+1} - 2[K^+]_e^{i} + [K^+]_e^{i-1}}{\Delta x^2} \, , \quad i = 2, \ldots , N-1 \, ,
    \end{equation}
    where $i$ indexes the array and $N$ denotes therefore its total number of elements (i.e. number of neurons and astrocytes) considered in the discrete network.

    Dirichlet boundaries conditions have been applied to the endpoints of the array, which aims to simulate healthy tissue. Thus, for the first and last indexes, we must impose such conditions to obtain valid approximations there. Indeed,
    \begin{equation}
        \begin{split}
            &\text{First index approximation: } \frac{\partial^2 [K^+]_e}{\partial x^2} \approx \frac{[K^+]_e^2 - 2 [K^+]_e^1 + [K^+]_{e,\text{bdry}}}{\Delta x^2} \,, \\[0.3em]
            &\text{Last index approximation: } \frac{\partial^2 [K^+]_e}{\partial x^2} \approx \frac{[K^+]_{e,\text{bdry}} - 2 [K^+]_e^N + [K^+]_e^{N-1}}{\Delta x^2} \, .
        \end{split}
    \end{equation}
    Same approximations apply to the diffusion term in the $[Na^+]_e$ equation for system \eqref{gHuguet_model}.

   The initial conditions are identical for all the pairs.
    We used the stiff integrator \texttt{ode15s} of \textsc{Matlab} with crude tolerances.

\section{Computation of the derivatives of the functions $H^{\ast}$ } \label{ap:derH}

    In this section we present the computations of the derivatives of the functions $H^{\ast}(z)$, for $\ast=l,m,r$ given in \eqref{eq:H*} and plotted in Figure~\ref{equilibrium_points_slow_system}. Thus,
    \begin{equation}
        \frac{dH^{\ast}}{dz}(z)=\frac{\partial h}{\partial x}(X^\ast(z),Y(z),z) \, \frac{d X^{\ast}}{d z}(z) + \frac{\partial h}{\partial y}(X^{\ast}(z),Y(z),z) \, \frac{d Y}{d z}(z) + \frac{\partial h}{\partial z}(X^{\ast}(z),Y(z),z) \,,
    \end{equation}
    where $h(x,y,z)$ is the function given in \eqref{eq:fgh}. The first partial derivatives of $h$ are
    \begin{equation} \label{derivatives_function_h}
        \begin{split}
            \frac{\partial h}{\partial x} &= \frac{10 S_N}{F \Omega_e} \big( g_K 4 n_{\infty}^3 n_{\infty}' (x - V_K) + g_K n_{\infty}^4 \big) \, , \\[0.2em]
            \frac{\partial h}{\partial y} &= \frac{10 S_A}{F \Omega_e} P_K \frac{F^2}{R T} \bigg( \frac{z e^{-\phi} - [K^+]_i^A}{e^{-\phi} - 1} + \phi e^{-\phi} \frac{z - [K^+]_i^A}{(e^{-\phi} - 1)^2} \bigg) \, , \\[0.2em]
            \frac{\partial h}{\partial z} &= 
            \frac{10 S_N}{F \Omega_e} \Bigg( -g_K n_{\infty}^4 \frac{R T}{F} \frac{[K^+]_i + \alpha_0 (1 + \Omega_a/\Omega_n) z}{z [K^+]_i} - 4 \rho_N \frac{z K_K}{(z + K_K)^3} \bigg( \frac{[Na^+]_i}{K_{Na} + [Na^+]_i} \bigg)^3 \Bigg) \\[0.2em] &\quad + \frac{10 S_A}{F \Omega_e} \Bigg( P_K F \frac{\phi e^{-\phi}}{e^{-\phi} - 1} - 4 \rho_A \frac{z K_K}{(z + K_K)^3} \bigg( \frac{[Na^+]_i^A}{K_{Na,A} + [Na^+]_i^A} \bigg)^3 \Bigg)\, ,
        \end{split}
    \end{equation}
    where $\phi = \frac{y F}{R T}$. For $\phi=0$ we use L'H\^opital's rule to make the derivatives continuous.
    
    The expressions for the derivatives of the functions $X^{\ast}(z)$ and $Y(z)$ are obtained from the differentiation of the equations $f(X^{\ast}(z),z) = 0$, for $\ast=l,m,r$, and $g(Y(z),z) = 0$ with respect to $z$. Indeed,
    \[\frac{d X^{\ast}}{d z} = -\dfrac{\nicefrac{\partial f}{\partial z}}{\nicefrac{\partial f}{\partial x}}, \quad \textrm{and} \quad  \frac{d Y}{d z} = -\frac{\nicefrac{\partial g}{\partial z}}{\nicefrac{\partial g}{\partial y}} \, .\]

\section{Taylor expansions of the manifold $\M_{\varepsilon}^r$ in $\varepsilon$} \label{sec:expansionmanifold}

To compute Taylor expansions of the manifold $\M_{\varepsilon}^r$ in $\varepsilon$ we solve the invariance equation \eqref{invariance_equation_xy} by matching powers of $\varepsilon$.

Developing functions $f$, $g$ and $h$ as Taylor series around the coefficients of order 0, $m_0$ and $n_0$, yields
    \begin{equation} \label{taylor_series_fgh}
        \begin{split}
            f(m,z) &= f(m_0,z) + \frac{\partial f}{\partial x}(m_0,z) (\varepsilon m_1 + \varepsilon^2 m_2) + \frac{1}{2}\frac{\partial^2 f}{\partial x^2}(m_0,z) \varepsilon^2 m_1^2 + \mathcal{O}(\varepsilon^3) \, , \\[0.2em]
            g(n,z) &= g(n_0,z) + \frac{\partial g}{\partial y}(n_0,z) (\varepsilon n_1 + \varepsilon^2 n_2) + \frac{1}{2}\frac{\partial^2 g}{\partial y^2}(n_0,z) \varepsilon^2 n_1^2 + \mathcal{O}(\varepsilon^3) \, , \\[0.2em]
            h(m,n,z) &= h(m_0,n_0,z) + \frac{\partial h}{\partial x}(m_0,n_0,z) \varepsilon m_1 + 
\frac{\partial h}{\partial y}(m_0,n_0,z) \varepsilon n_1  + \mathcal{O}(\varepsilon^2) \, .
        \end{split}
    \end{equation}

    After substituting the Taylor expansions \eqref{taylor_series_fgh} into system \eqref{invariance_equation_xy}, we are able to solve the resulting system of equations up to any power of $\varepsilon$:
    
\begin{itemize}
        \item $\mathcal{O}(1)$ :
        \begin{equation*}
            \begin{split}
                f(m_0,z) &= 0 \quad \Longrightarrow \quad m_0(z,w) = X^r(z) \, , \,  \quad \text{(see red curve in Figure \ref{critical_manifold_dissection})} \\[0.2em]
                g(n_0,z) &= 0 \quad \Longrightarrow \quad n_0(z,w) = Y(z) \, . \quad \text{(see orange curve in Figure \ref{critical_manifold_dissection})}
            \end{split}
        \end{equation*}
        \item $\mathcal{O}(\varepsilon)$ :
        \begin{equation*}
            \begin{split}
                \frac{\partial f}{\partial x}(m_0,z) \,m_1 &= c w \frac{\partial m_0}{\partial z} + \cancelto{0}{c^2 w \frac{\partial m_0}{\partial w}} - \cancelto{0}{c \, h(m_0,n_0,z) \frac{\partial m_0}{\partial w}} \quad \Longrightarrow \quad m_1(z,w) =  -\frac{c w \, f_z(X^r(z),z)}{\big(f_x(X^r(z),z)\big)^2} \, , \\[0.2em]
                \frac{\partial g}{\partial y}(n_0,z) \, n_1 &= c w \frac{\partial n_0}{\partial z} + \cancelto{0}{c^2 w \frac{\partial n_0}{\partial w}} - \cancelto{0}{c \, h(m_0,n_0,z) \frac{\partial n_0}{\partial w}} \quad \Longrightarrow \quad n_1(z,w) = -\frac{c w \, g_z(Y(z),z)}{\big(g_y(Y(z),z)\big)^2} \, .
            \end{split}
        \end{equation*}

        Observe that, since $f(X^r(z),z)=0$ and $g(Y(z),z)=0$, we have used
    \begin{equation}\label{eq:pmn0z}
        \frac{\partial m_0}{\partial z} = \frac{d X^r(z)}{dz} = -\frac{f_z(m_0,z)}{f_x(m_0,z)} \quad \text{and} \quad \frac{\partial n_0}{\partial z} = \frac{d Y}{dz}= -\frac{g_z(n_0,z)}{g_y(n_0,z)} \, .
    \end{equation}

        \item $\mathcal{O}(\varepsilon^2)$ :
        \begin{equation*}
            \begin{split}
                &\begin{aligned}
                    f_x(m_0,z) \, m_2 + \frac{1}{2!} f_{x x}(m_0,z) \, m_1^2 &= c w \frac{\partial m_1}{\partial z} + c^2 w \frac{\partial m_1}{\partial w} - c h(m_0,n_0,z) \frac{\partial m_1}{\partial w} \\[0.2em] & \quad - c \big( h_x(m_0,n_0,z) \, m_1 + h_y(m_0,n_0,z) \, n_1 \big) \cancelto{0}{\frac{\partial m_0}{\partial w}} \, .
                \end{aligned}
                \\[1em]
                & \Longrightarrow \quad m_2(z,w) = -\frac{1}{2!} \frac{f_{x x}(m_0,z)}{f_x(m_0,z)} \, m_1^2 + \frac{c w}{f_x(m_0,z)} \frac{\partial m_1}{\partial z} - c \big(c w - h(m_0,n_0,z) \big) \frac{c \, f_z(m_0,z)}{\big(f_x(m_0,z)\big)^3} \, ,
                \\[1em]
                &\begin{aligned}
                    g_y(n_0,z) \, n_2 + \frac{1}{2!} g_{y y}(n_0,z) \, n_1^2 &= c w \frac{\partial n_1}{\partial z} + c^2 w \frac{\partial n_1}{\partial w} - c h(m_0,n_0,z) \frac{\partial n_1}{\partial w} \\[0.2em] & \quad - c \big( h_x(m_0,n_0,z) \, m_1 + h_y(m_0,n_0,z) \, n_1 \big) \cancelto{0}{\frac{\partial n_0}{\partial w}} \, .
                \end{aligned}
                \\[1em]
                & \Longrightarrow \quad n_2(z,w) = -\frac{1}{2!} \frac{g_{y y}(n_0,z)}{g_y(n_0,z)} \, n_1^2 + \frac{c w}{g_y(n_0,z)} \frac{\partial n_1}{\partial z} - c \big(c w - h(m_0,n_0,z) \big) \frac{c \, g_z(n_0,z)}{\big(g_y(n_0,z)\big)^3} \, .
            \end{split}
        \end{equation*}
    \end{itemize}
    where
    \begin{equation}\label{eq:pmn1z}
        \begin{split}
            \frac{\partial m_1}{\partial z} &= -c w \frac{(f_{xz}(m_0,z) \frac{\partial m_0}{\partial z} + f_{zz}(m_0,z)) \, f_x(m_0,z) - 2 f_z(m_0,z) (f_{xx}(m_0,z) \frac{\partial m_0}{\partial z} + f_{zx}(m_0,z))}{(f_x(m_0,z))^3} \, , \\[0.2em]
            \frac{\partial n_1}{\partial z} &= -c w \frac{(g_{yz}(n_0,z) \frac{\partial n_0}{\partial z} + g_{zz}(n_0,z)) \, g_y(n_0,z) - 2 g_z(n_0,z) (g_{yy}(n_0,z) \frac{\partial n_0}{\partial z} + g_{zy}(n_0,z))}{(g_y(n_0,z))^3} \, ,
        \end{split}
    \end{equation}

\section{Computation of the partial derivatives of the functions $\widetilde f(x,z)$ and $\widetilde g(y,z)$ }

    In this section we detail the expression of the first partial derivatives of the functions $\widetilde f(x,z)$ and $\widetilde g(y,z)$ given in \eqref{eq:fgh}.
    The derivatives of higher order can be computed analogously.
    Recall that $x=V_N$, $y=V_A$ and $z=[K^+]_e$. The term $1/C_m$ has been omitted because, in practice, $C_m$ is equal to 1.
    Thus,
   \begin{equation} \label{derivatives_functions_f_g}
        \begin{split}
            &\begin{aligned}
                \frac{\partial \widetilde f}{\partial x} =& \, - \bigg[ g_{Na} \, 3 \, m_\infty^2 \, m_\infty' \, (1 - n_\infty) \, (x - V_{Na}) + g_{Na} \, m_\infty^3 \, (1 - n_\infty' x - n_\infty + n_\infty' V_{Na}) \\[0.2em] &\quad+ g_{NaP} \, m_\infty' \, h_p \, (x - V_{Na}) + g_{NaP} \, m_{p_\infty} \, h_p \\[0.2em] &\quad+ g_K \, 4 \, n_\infty^3 \, n_\infty' \, (x - V_K) + g_K \, n_\infty^4 + g_L \bigg] \, ,
            \end{aligned}
            \\[0.5em]
            &\begin{aligned}
                \frac{\partial \widetilde f}{\partial z} =& \, - \bigg[ -g_{K} \, n_\infty^4 \frac{R T}{F} \frac{[K^+]_i + \alpha_0 (1 + \Omega_a/\Omega_n) z}{z [K^+]_i} + 2 \, \rho_N \frac{z K_K}{(z + K_K)^3} \bigg( \frac{[Na^+]_i}{K_{Na} + [Na^+]_i} \bigg)^3 \bigg] \, ,
            \end{aligned}
            \\[0.5em]
            &\begin{aligned}
                \frac{\partial \widetilde g}{\partial y} =& \, - \bigg[ P_{Na} \frac{F^2}{RT} \bigg( \frac{[Na^+]_e e^{-\phi} - [Na^+]_i^A}{e^{-\phi} - 1} + \phi e^{-\phi} \frac{[Na^+]_e - [Na^+]_i^A}{(e^{-\phi} - 1)^2} \bigg) \\[0.2em]
                &\quad+ P_{K} \frac{F^2}{RT} \bigg( \frac{z e^{-\phi} - [K^+]_i^A}{e^{-\phi} - 1} + \phi e^{-\phi} \frac{z - [K^+]_i^A}{(e^{-\phi} - 1)^2} \bigg) \bigg] \, ,
            \end{aligned}
            \\[0.5em]
            &\begin{aligned}
                \frac{\partial \widetilde g}{\partial z} =& \, - \bigg[ 2 \rho_A \frac{z K_{K,A}}{(z + K_{K,A})^3} \bigg( \frac{[Na^+]_i^A}{K_{Na,A} + [Na^+]_i^A} \bigg)^3 + P_K F \frac{\phi e^{-\phi}}{e^{-\phi} - 1} \bigg] \, ,
            \end{aligned}
        \end{split}
    \end{equation}
    where $\phi := \frac{y F}{RT}$. 

    We also define the derivatives when $y=0$ using L'H\^opital's rule, and we obtain
    \[\frac{\partial \widetilde g}{\partial y} (y,z)\bigg|_{y=0}= -\frac{1}{2} \frac{F^2}{R T} \Big(P_{Na} \big([Na^+]_e + [Na^+]_i^A\big) + P_{K} \big(z + [K^+]_i^A\big) \Big),\]
    and
    \[
    \frac{\partial \widetilde g}{\partial z}(y,z)\bigg|_{y=0}=- \bigg[ 2 \rho_A \frac{z K_{K,A}}{(z + K_{K,A})^3} \bigg( \frac{[Na^+]_i^A}{K_{Na,A} + [Na^+]_i^A} \bigg)^3 - P_K F \bigg] \, .
    \]

\section{Linear local approximation for the 1D stable manifold \texorpdfstring{$\widetilde \Gamma^s$}{}}\label{ap:DF} 

To compute the eigenvalues and first order approximation of the 1D stable manifold $\widetilde \Gamma^s$ we need to compute the differential matrix for the restricted system \eqref{slow_subsystem_like}, referred as $F_2$, given by
    \begin{equation}
        DF_2(z,w) = 
        \begin{pmatrix}
            0 & \varepsilon \\[0.5em]
            -\varepsilon \Big(\dfrac{\partial h}{\partial x} \dfrac{\partial x}{\partial z} + \dfrac{\partial h}{\partial y} \dfrac{\partial y}{\partial z} + \dfrac{\partial h}{\partial z} \Big) & \varepsilon \bigg(c - \Big(\dfrac{\partial h}{\partial x} \dfrac{\partial x}{\partial w} + \dfrac{\partial h}{\partial y} \dfrac{\partial y}{\partial w} \Big) \bigg)
        \end{pmatrix} \, ,
    \end{equation}
    where $x$ and $y$ are defined as functions of $z$ and $w$ in \eqref{eq:embed}. Thus, 
    \begin{equation*}
        \begin{split}
            \frac{\partial x}{\partial z} &= \frac{\partial}{\partial z} (m_0 + \varepsilon m_1 + \varepsilon^2 m_2) = -\frac{f_z(m_0,z)}{f_x(m_0,z)} + \varepsilon \frac{\partial m_1}{\partial z} + \varepsilon^2 \frac{\partial m_2}{\partial z} \, , \\[0.2em]
            \frac{\partial x}{\partial w} &= \frac{\partial}{\partial w} (m_0 + \varepsilon m_1 + \varepsilon^2 m_2) = -\varepsilon \frac{c f_z(m_0,z)}{(f_x(m_0,z))^2} + \varepsilon^2 \frac{\partial m_2}{\partial w} \, , \\[0.2em]
            \frac{\partial y}{\partial z} &= \frac{\partial}{\partial z} (n_0 + \varepsilon n_1 + \varepsilon^2 n_2) = -\frac{g_z(n_0,z)}{g_y(n_0,z)} + \varepsilon \frac{\partial n_1}{\partial z} + \varepsilon^2 \frac{\partial n_2}{\partial z} \, , \\[0.2em]
            \frac{\partial y}{\partial w} &= \frac{\partial}{\partial w} (n_0 + \varepsilon n_1 + \varepsilon^2 n_2) = -\varepsilon \frac{c g_z(n_0,z)}{(g_y(n_0,z))^2} + \varepsilon^2 \frac{\partial n_2}{\partial w} \, ,
        \end{split}
    \end{equation*}
    where the derivatives $\partial m_0/\partial z$ and $\partial n_0/\partial z$ are taken from \eqref{eq:pmn0z} and we have used $\partial m_1/\partial w=\partial n_1/\partial w=0$. Finally, $\partial m_1/\partial z$ and $\partial n_1/\partial z$ are given in \eqref{eq:pmn1z}. For the other ones, we have
    \begin{equation*}
        \begin{split}
        &\begin{aligned}
            \frac{\partial m_2}{\partial z} &= -\frac{1}{2} \bigg[\frac{\big( f_{xxx}(m_0,z) \partial_z m_0 + f_{zxx}(m_0,z) \big) f_x(m_0,z) - f_{xx}(m_0,z)\big(f_{xx}(m_0,z) \partial_z m_0 + f_{zx}(m_0,z)\big)}{(f_x(m_0,z)^2} m_1^2 \\[0.2em]
            &\quad + \frac{f_{xx}(m_0,z)}{f_x(m_0,z)} 2 m_1 \frac{\partial m_1}{\partial z}\bigg] + \bigg[\frac{-c w \big(f_{xx}(m_0,z) \partial_z m_0 + f_{zx}(m_0,z)\big)}{(f_x(m_0,z))^2} \frac{\partial m_1}{\partial z} + \frac{c w}{f_x(m_0,z)} \frac{\partial^2 m_1}{\partial z^2}\bigg] \\[0.2em]
            &\quad - c \bigg[\big(-h_x(m_0,n_0,z) \partial_z m_0 - h_y(m_0,n_0,z) \partial_z n_0 - h_z(m_0,n_0,z)\big)\frac{c f_z(m_0,z)}{(f_x(m_0,z))^3} \\[0.2em] &\quad + (c w - h(m_0,n_0,z)) c \frac{\big(f_{xz}(m_0,z) \partial_z m_0 + f_{zz}(m_0,z)\big) f_x(m_0,z) - 3 f_z(m_0,z) \big(f_{xx}(m_0,z) \partial_z m_0 + f_{zx}(m_0,z)\big)}{(f_x(m_0,z))^4}\bigg] \, ,
        \end{aligned} \\[0.3em]
        &\frac{\partial m_2}{\partial w} = -\frac{1}{2} \frac{f_{xx}(m_0,z)}{f_x(m_0,z)} 2 m_1 \bigg(-c \frac{f_z(m_0,z)}{(f_x(m_0,z))^2}\bigg) + \frac{2 c}{f_x(m_0,z)} \frac{\partial m_1}{\partial z} - c^3 \frac{f_z(m_0,z)}{(f_x(m_0,z))^3} \, ,
        \end{split}
    \end{equation*}
    with
    \begin{equation*}
        \begin{split}
            \frac{\partial^2 m_1}{\partial z^2} &= (-cw) \frac{1}{(f_x(m_0,z))^4} \Biggl\{\bigg[ \bigg( \big(f_{xxz}(m_0,z) \partial_z m_0 + f_{zxz}(m_0,z)\big) \partial_z m_0 \\[0.2em]
            &\quad + f_{xz}(m_0,z) \partial^2_{zz} m_0 + \big( f_{xzz}(m_0,z) \partial_z m_0 + f_{zzz}(m_0,z) \big) \bigg) f_x(m_0,z) \\[0.2em]
            &\quad - \big(f_{xz}(m_0,z) \partial_z m_0 + f_{zz}(m_0,z)\big)\big(f_{xx}(m_0,z) \partial_z m_0 + f_{zx}(m_0,z)\big) \\[0.2em] &\quad - 2 f_z(m_0,z) \bigg( \big( f_{xxx}(m_0,z) \partial_z m_0 + f_{zxx}(m_0,z) \big) \partial_z m_0 + f_{xx}(m_0,z) \partial_{zz}^2 m_0 \\[0.2em] &\quad + \big(f_{xzx}(m_0,z) \partial_z m_0 + f_{zzx}(m_0,z) \big) \bigg) \bigg] f_x(m_0,z) \\[0.2em] &\quad - 3 \big(f_{xx}(m_0,z) \partial_z m_0 + f_{zx}(m_0,z)\big) \bigg( \big(f_{xz}(m_0,z) \partial_z\ m_0 + f_{zz}(m_0,z) \big) f_x(m_0,z) \\[0.2em] &\quad - 2 f_z(m_0,z) \big( f_{xx}(m_0,z) \partial_z m_0 + f_{zx}(m_0,z) \big) \bigg) \Biggr\} \, .
        \end{split} 
    \end{equation*}
    One obtains the expressions for the partial derivatives of $n_2$ with respect to $z$ and $w$ by carefully exchanging $f$'s and $x$'s for $g$'s and $y$'s, respectively.
    
\section*{Acknowledgments}
Work produced with the support of the grant PID-2021-122954NB-100 funded by MCIN/AEI/ 10.13039/501100011033 and “ERDF: A way of making Europe”. T.M.S and G.H acknowledge the Maria de Maeztu Award for Centers and Units of Excellence in R\&D (CEX2020-001084-M). 
T. M. S. is supported by the Catalan Institution for Research and Advanced Studies via an ICREA Academia Prize 2019. 
We also acknowledge the use of the UPC Dynamical Systems group’s cluster for research computing \texttt{https://dynamicalsystems.upc.edu/en/computing/}

\end{document}